\documentclass[12pt]{article}

\usepackage{amsmath}
\usepackage{amsfonts}
\usepackage{latexsym}
\usepackage{graphicx}

\newtheorem{proposition}{Proposition}[section]

\newtheorem{corollary}{Corollary}[section]
\newtheorem{remark}{Remark}[section]
\newtheorem{example}{Example}[section]

\date{ }

\begin{document}

\title{Topological optimization via cost penalization}

\author{Cornel Marius Murea$^1$, Dan Tiba$^2$\\
{\normalsize $^1$ D\'epartement de Math\'ematiques, IRIMAS,}\\
{\normalsize Universit\'e de Haute Alsace, France,}\\
{\normalsize cornel.murea@uha.fr}\\
{\normalsize $^2$ Institute of Mathematics (Romanian Academy) and}\\ 
{\normalsize Academy of Romanian Scientists, Bucharest, Romania,}\\ 
{\normalsize dan.tiba@imar.ro}
}

\maketitle

\begin{abstract}%
We consider general shape optimization problems governed by Dirichlet boundary value problems.
The proposed approach may be extended to other boundary conditions as well. It is based on a recent 
representation result for implicitly defined manifolds, due to the authors, and it is formulated 
as an optimal control problem. The discretized approximating problem is introduced and we give 
an explicit construction of the associated discrete gradient. Some numerical examples are also indicated.

\textbf{Keywords: geometric optimization; optimal design; topological variations; optimal control methods; 
discrete gradient}
\end{abstract}

\section{Introduction}
\setcounter{equation}{0}

Shape optimization is a relatively young branch of mathematics, with important
modern applications in engineering and design. Certain optimization problems in 
mechanics, thickness optimization for plate or rods, geometric optimization
of shells, curved rods, drag minimization in fluid mechanics, etc are some examples.
Many appear naturally in the form of control by coefficients problems, due
to the formulation of the mechanical models, with the geometric characteristics
entering the coefficients of the differential operators. See \cite{Tiba2006}, Ch. 6,
where such questions are discussed in details.

It is the aim of this article to develop an optimal control approach, using
penalization methods, to general shape optimization problems as investigated in 
\cite{Pironneau1984}, \cite{Sokolowski1992}, \cite{Bucur2005}, \cite{Henrot2005}, 
\cite{Delfour2001}, etc. We underline that our methodology allows simultaneous topological 
and boundary variations.

Here, we fix our attention on the case of Dirichlet boundary conditions and we study 
the typical problem (denoted by ($\mathcal{P}$)):
\begin{eqnarray}
\min_\Omega \int_E j\left(\mathbf{x}, y_\Omega(\mathbf{x})\right)d\mathbf{x},
\label{1.1}\\
-\Delta y_\Omega = f \hbox{ in }\Omega,
\label{1.2}\\
y_\Omega = 0  \hbox{ on }\partial\Omega,
\label{1.3}
\end{eqnarray}
where $E\subset\subset D \subset \mathbb{R}^2$ are given bounded domains,
$D$ is of class $\mathcal{C}^{1,1}$ and the minimization parameter,
the unknown domain $\Omega$, satisfies $E \subset \Omega \subset D$
and other possible conditions defining a class of admissible domains.
Notice that the case of dimension 
two is of interest in shape optimization.
Moreover, $f\in L^2(D)$, $j:D\times\mathbb{R}\rightarrow\mathbb{R}$ is
some Carath\'eodory mapping. More assumptions
or constraints will be imposed later. Other boundary conditions or differential operators
may be handled as well via this control approach and we shall examine such questions in 
a subsequent paper.

For fundamental properties and methods in optimal control theory, we quote
\cite{Lions1969}, \cite{Clarke1983}, \cite{Barbu1986}, \cite{Tiba1994}.
The problem (\ref{1.1})-(\ref{1.3}) and its approximation are strongly non convex and
challenging both from the numerical and theoretical points of view.
The investigation from this paper continues the one
in \cite{Tiba2018a} and is essentially based on the recent implicit 
parametrization method as developed in \cite{Tiba2018}, \cite{Tiba2015}, \cite{Tiba2013}, 
that provides an efficient analytic representation of the unknown domains.

The Hamiltonian approach to implicitly defined manifolds will be briefly recalled
together with other preliminaries in Section 2.
The precise formulation of the problem and its approximation is analyzed in Section 3
together with its differentiability properties. In Section 4, we study the discretized version 
and find the general form of the discrete gradient.
The last section is devoted to some numerical experiments, using this paper approach.

The method studied in this paper has a certain complexity due to the use of Hamiltonian systems and its main advantage is the possibility to extend it to other boundary conditions or boundary observation problems. This will be performed in a subsequent article.

\section{Preliminaries}
\setcounter{equation}{0}

Consider the Hamiltonian system
\begin{eqnarray}
x_1^\prime(t) & = & -\frac{\partial g}{\partial x_2}\left(x_1(t),x_2(t)\right),\quad t\in I,
\label{2.1}\\
x_2^\prime(t) & = &  \frac{\partial g}{\partial x_1}\left(x_1(t),x_2(t)\right),\quad t\in I,
\label{2.2}\\
\left(x_1(0),x_2(0)\right)& = &\left(x_1^0,x_2^0\right),
\label{2.3}
\end{eqnarray}
where $g:D\rightarrow\mathbb{R}$ is in $\mathcal{C}^1(\overline{D})$,
$\left(x_1^0,x_2^0 \right)\in D$ and $I$ is the local existence interval
for (\ref{2.1})--(\ref{2.3}), around the origin, obtained via the Peano theorem.
The conservation property \cite{Tiba2013} of the Hamiltonian gives:

\begin{proposition}\label{prop:2.1}
We have
\begin{equation}\label{2.4}
g\left(x_1(t),x_2(t)\right)=g\left(x_1^0,x_2^0\right), \quad t\in I.
\end{equation}  
\end{proposition}

In the sequel, we assume that
\begin{equation}\label{2.5}
g\left(x_1^0,x_2^0\right)=0,\quad \nabla g\left(x_1^0,x_2^0\right)\neq 0.
\end{equation}
Under condition (\ref{2.5}), i.e. in the noncritical case, it is known that the solution
of (\ref{2.1})--(\ref{2.3}) is also unique (by applying the implicit functions
theorem to (\ref{2.4}), \cite{4}).

\begin{remark}\label{rem:1}
  In higher dimension, iterated Hamiltonian systems were introduced in \cite{Tiba2018}
  and uniqueness and regularity properties are proved. Some relevant examples
  in dimension three are discussed in \cite{Tiba2015}. In the critical case, generalized solutions can be obtained \cite{Tiba2013}, \cite{Tiba2018}.
\end{remark}

We define now the family $\mathcal{F}$ of admissible functions
$g\in \mathcal{C}^2(\overline{D})$ that satisfy the conditions:
\begin{eqnarray}
g(x_1,x_2) & > & 0,\quad\hbox{on }\partial D,
\label{2.6}\\
|\nabla g(x_1,x_2)| & > &  0,\quad\hbox{on }
\mathcal{G}=\left\{ (x_1,x_2)\in D;\ g(x_1,x_2)=0\right\},
\label{2.7}\\
g(x_1,x_2) & < & 0,\quad\hbox{on } \overline{E}.
\label{2.8}
\end{eqnarray}

Condition (\ref{2.6}) says that $\mathcal{G}\cap \partial D=\emptyset$ and condition
(\ref{2.7}) is an extension of (\ref{2.5}). In fact, it is related to the hypothesis
on the non existence of equilibrium points in the Poincare-Bendixson theorem, \cite{hsd},
Ch. 10, and the same is valid for the next proposition.
The family $\mathcal{F}$ defined by (\ref{2.6})--(\ref{2.8}) is obviously very rich, but
it is not ``closed'' 
(we have strict inequalities). Our approach here, gives a descent algorithm for the shape
optimization problem 
($\mathcal{P}$) and existence of optimal shapes is not discussed.

Following \cite{Tiba2018a}, we have the following two propositions:

\begin{proposition}\label{prop:2.2}
  Under hypotheses (\ref{2.6}), (\ref{2.7}), $\mathcal{G}$ is a finite union of closed curves of
  class $\mathcal{C}^2$, without self intersections, parametrized by
  (\ref{2.1})--(\ref{2.3}), when some initial point $(x_1^0,x_2^0)$ is chosen on each
  component of $\mathcal{G}$.
\end{proposition}

If $r\in \mathcal{F}$ as well, we define the perturbed set
\begin{equation}\label{2.9}
\mathcal{G}_\lambda=\left\{ (x_1,x_2)\in D;\ (g+\lambda r)(x_1,x_2)=0,\ \lambda\in \mathbb{R}\right\}.
\end{equation}
We also introduce the neighborhood $V_\epsilon$, $\epsilon>0$
\begin{equation}\label{2.10}
V_\epsilon=\left\{ (x_1,x_2)\in D;\ d[(x_1,x_2),\mathcal{G}] < \epsilon\right\},
\end{equation}
where $d[(x_1,x_2),\mathcal{G}]$ is the distance from a point to $\mathcal{G}$.

\begin{proposition}\label{prop:2.3}
  If $\epsilon>0$ is small enough, there is $\lambda(\epsilon)>0$ such that,
  for $|\lambda| < \lambda(\epsilon)$ we have $\mathcal{G}_\lambda \subset V_\epsilon$
  and $\mathcal{G}_\lambda$ is a finite union of class $\mathcal{C}^2$ closed curves.
\end{proposition}

\begin{remark}\label{rem:2}
  The inclusion $\mathcal{G}_\lambda \subset V_\epsilon$ shows that $\mathcal{G}_\lambda\rightarrow \mathcal{G}$ for
  $\lambda\rightarrow 0$, in the Hausdorff-Pompeiu metric \cite{Tiba2006}.
  In a "small" neighborhood of each component of $\mathcal{G}$ there is exactly one component
  of $\mathcal{G}_\lambda$ if $|\lambda| < \lambda (\epsilon)$, due to this convergence
  property and the implicit functions theorem applied in the initial condition of the
  perturbed Hamiltonian system derived from (\ref{2.1})--(\ref{2.3}) .
\end{remark}

\begin{proposition}\label{prop:2.4}
Denote by $T_g$, $T_{g+\lambda r}$ the periods of the trajectories of (\ref{2.1})--(\ref{2.3}), 
corresponding to $g$, $g+\lambda r$ respectively. Then $T_{g+\lambda r} \rightarrow T_g$ as
$\lambda \rightarrow 0$.
\end{proposition}

\noindent
\textbf{Proof.}
If $(x_1,x_2)$, respectively $(x_{1\lambda},x_{2\lambda})$ are the corresponding trajectories of
(\ref{2.1})--(\ref{2.3}) respectively, then they are bounded by Proposition \ref{prop:2.3},
if $|\lambda| < \lambda(\epsilon)$.
Consequently, $\nabla g$ may be assumed Lipschitzian with constant $L_g$ and we have
\begin{equation}\label{2.11}
|(x_1,x_2)-(x_{1\lambda},x_{2\lambda}) |(t) \leq
\lambda C + L_g \int_0^t |(x_1,x_2)-(x_{1\lambda},x_{2\lambda}) | dt,
\end{equation}
where  we also use that $\nabla r(x_{1\lambda},x_{2\lambda})$ is bounded since 
$(x_{1\lambda},x_{2\lambda})$ is bounded on $\mathbb{R}$.
We infer by (\ref{2.11}) that
\begin{eqnarray}
|(x_1,x_2)-(x_{1\lambda},x_{2\lambda}) |(t) \leq\lambda\, ct, t\in \mathbb{R},\label{2.12}\\
|(x_1^\prime,x_2^\prime)-(x_{1\lambda}^\prime,x_{2\lambda}^\prime) |(t) \leq\lambda\, ct, t\in \mathbb{R},\label{2.13}
\end{eqnarray}
for $|\lambda| < \lambda(\epsilon)$ and with some constant independent of $\lambda$, by Gronwall lemma.

Both trajectories start from $(x_1^0,x_2^0)$, surround $\overline{E}$, have no self intersections ( but $(x_{1\lambda},x_{2\lambda})$ may intersect $(x_1,x_2)$ even on infinity of times).
We study them on $[0,jT_g]$, $j<2$, for instance.

Assume that $(x_{1\lambda},x_{2\lambda})$ has the period $T_{g+\lambda r} > jT_g$. Since $(x_1,x_2)$ is periodic with period $T_g$ and
relations (\ref{2.12})--(\ref{2.13}) show that $(x_{1\lambda},x_{2\lambda})$ is very close to $(x_1,x_2)$ in every $t \in [0,jT_g]$ it yields that $(x_{1\lambda},x_{2\lambda})$ is, as well, surrounding $\overline{E}$ at least once.
As it may have no self intersections, it yields that $(x_{1\lambda},x_{2\lambda})$ is as a limit cycle around $\overline{E}$.
Such arguments appear in the proof of the Poincar\'e-Bendixson
theorem, \cite{hsd}, Ch. 10.
That is $(x_{1\lambda},x_{2\lambda})$ cannot be periodic - and
this is a false conclusion due to Proposition \ref{prop:2.3}.

Consequently, we get
\begin{equation}\label{2.14}
T_{g+\lambda r} \leq jT_g,\ |\lambda| < \lambda(\epsilon).
\end{equation}
On a subsequence, by (\ref{2.14}), we obtain $T_{g+\lambda r} \rightarrow T^* \leq jT_g$. We assume that $T^* \neq T_g$.
It is clear that $\left(x_{1}(T^*),x_{2}(T^*)\right)\neq
\left(x_{1}(T_g),x_{2}(T_g)\right)$ by the definition of the
period. However, relation (\ref{2.12}) and the related
convergence properties give the opposite conclusion.
This contradiction shows that $T^* = T_g$ and the convergence
is valid on the whole sequence.\quad$\Box$

\begin{remark}\label{rem:2.7}
  Usually, the perturbation of a periodic solution may not
  be periodic and just asymptotic convergence properties are
  valid, under certain assumptions, Sideris
  \cite{Sideris2013}.
  A natural
  question, taking into account (\ref{2.12})--(\ref{2.13}),
  is whether $|T_{g+\lambda r}-T_g| \leq c |\lambda|$, for
  $c>0$ independent of $\lambda$,
  $|\lambda| < \lambda(\epsilon)$.
\end{remark}

\section{The optimization problem and its approximation}
\setcounter{equation}{0}

Starting from the family $\mathcal{F}$ of admissible functions, we define
the family $\mathcal{O}$ of admissible domains in the shape optimization
problem  (\ref{1.1})-(\ref{1.3}) as the connected component of the open set
$\Omega_g$, $g\in \mathcal{F}$
\begin{equation}\label{3.1}
\Omega_g=\left\{ (x_1,x_2)\in D;\ g(x_1,x_2) < 0\right\}
\end{equation}
that contains $E$. Clearly $E\subset \Omega_g$ by (\ref{2.8}).
Notice as well that the domain $\Omega_g$ defined by (\ref{3.1})
(we use this notation for the domain as well) is not simply connected, in general.
This is the reason why the approach to (\ref{1.1})-(\ref{1.3}) that we discuss here is related to
topological optimization in optimal design problems.
But, it also combines topological and boundary variations.

The penalized problem, $\epsilon >0$, is given by:
\begin{equation}\label{3.2}
\min_{g\in \mathcal{F},\, u\in L^2(D)} 
\left\{ 
\int_E j\left(\mathbf{x}, y_\epsilon(\mathbf{x})\right)d\mathbf{x}
+\frac{1}{\epsilon}\int_{I_g} \left(y_\epsilon(\mathbf{z}_g(t))\right)^2|\mathbf{z}_g^\prime(t)|dt
\right\}
\end{equation}
subject to
\begin{eqnarray}
-\Delta y_\epsilon & = & f+ (g+\epsilon)_+^2u,\quad\hbox{in } D,
\label{3.3}\\
y_\epsilon & = & 0,\quad\hbox{on } \partial D,
\label{3.4}
\end{eqnarray}
where $\mathbf{z}_g=(z_g^1,z_g^2)$ satisfies the Hamiltonian system (\ref{2.1})-(\ref{2.3})
in $I_g$ with some $(x_1^0,x_2^0)\in D\setminus\overline{E}$ such that $g(x_1^0,x_2^0)=0$
\begin{eqnarray}
\left(z_g^1\right)^\prime(t) & = & -\frac{\partial g}{\partial x_2}\left(\mathbf{z}_g(t)\right),\quad t\in I_g,
\label{3.5}\\
\left(z_g^2\right)^\prime(t) & = &  \frac{\partial g}{\partial x_1}\left(\mathbf{z}_g(t)\right),\quad t\in I_g,
\label{3.6}\\
\mathbf{z}_g(0)& = &\left(x_1^0,x_2^0\right)
\label{3.7}
\end{eqnarray}
and $I_g=[0,T_g]$ is the period interval for (\ref{3.5})-(\ref{3.7}), due to Proposition \ref{prop:2.2}.

The problem (\ref{3.2})-(\ref{3.7}) is an optimal control problem with controls $g\in \mathcal{F}$
and $u\in L^2(D)$ distributed in $D$. The state is given by 
$[y_\epsilon,z_g^1,z_g^2]\in H^2(D)\times \left(C^2(I_g)\right)^2$. We also have $y_\epsilon \in  H_0^1(D)$.
In case the corresponding domain $\Omega_g$ is not simply connected, in (\ref{3.7}) one has
to choose initial conditions on each component of $\partial\Omega_g$ and the penalization term becomes
 a finite sum due to Proposition \ref{prop:2.2}.
The method enters the class of fixed domain methods in shape optimization and can be compared with
\cite{Tiba1992}, \cite{Tiba2009}, \cite{Tiba2012}. It is essentially different from the
level set method
of Osher and Sethian \cite{Osher1988}, Allaire \cite{Allaire2002} or the SIMP approach of 
Bendsoe and Sigmund \cite{Bendsoe2003}.
From the computational point of view, it is easy to find initial condition (\ref{3.7}) on
each component
of $\mathcal{G}$ and the corresponding period intervals $I_g$ associated to (\ref{3.5})-(\ref{3.7}). See the last section as well.

We have the following general suboptimality property:

\begin{proposition}\label{prop:3.1}
Let $j(\cdot,\cdot)$ be a Carath\'eodory function on $D\times\mathbb{R}$, bounded from below
by a constant. Denote by $[y_n^\epsilon, g_n^\epsilon, u_n^\epsilon]$ a minimizing sequence in
(\ref{3.2})-(\ref{3.7}). Then, on a subsequence $n(m)$ the (not necessaryly admissible pairs)
$[\Omega_{g_{n(m)}^\epsilon}, y_{n(m)}^\epsilon]$ give a minimizing sequence in (\ref{1.1}),
$y_{n(m)}^\epsilon$ satisfies (\ref{1.2}) in
$\left\{ (x_1,x_2)\in D;\ g(x_1,x_2) < -\epsilon\right\}$ and (\ref{1.3}) is fulfilled with
a perturbation of order $\epsilon^{1/2}$.
\end{proposition}

\noindent
\textbf{Proof.} 
Let $[y_{g_m}, g_m]\in H^2(\Omega_{g_{m}})\times \mathcal{F}$ be a minimizing sequence in the problem
(\ref{1.1})-(\ref{1.3}), (\ref{3.1}). By the trace theorem, since $\partial \Omega_{g_{m}}$ and
$D$ are at least $\mathcal{C}^{1,1}$ under our assumptions, there is
$\widetilde{y}_{g_{m}}\in H^2(D)\cap H_0^1(D)$, not unique, such that
$\widetilde{y}_{g_{m}}=y_{g_m}$ in $\Omega_{g_{m}}$. We define the control $u_{g_{m}}\in L^2(D)$ as following:
\begin{eqnarray*}
u_{g_{m}} & = & 0,\quad\hbox{in } \Omega_{g_{m}},
\\
u_{g_{m}} & = & -\frac{\Delta \widetilde{y}_{g_{m}} + f}{(g_m+\epsilon)_+^2},
\quad\hbox{in } \partial D\setminus \Omega_{g_{m}},
\end{eqnarray*}
where $\Omega_{g_{m}}$ is the open set defined in (\ref{3.1}).
Notice that on the second line in the above formula, we have no singularity.
It is clear that the triple $[\widetilde{y}_{g_{m}},g_m,u_{g_{m}}]$ is admissible for the
problem (\ref{3.2})-(\ref{3.7}) with the same cost as in the original problem
(\ref{1.1})-(\ref{1.3}) since the penalization term in (\ref{3.2}) is null due
to the boundary condition (\ref{1.3}) satisfied by $\widetilde{y}_{g_{m}}$.
Consequently, there is $n(m)$ sufficient big, such that
\begin{eqnarray}
&&\int_E j\left(\mathbf{x}, y_{n(m)}^\epsilon(\mathbf{x})\right)d\mathbf{x}
+\frac{1}{\epsilon}\int_{I_{g_{n(m)}^\epsilon}} \left(y_{n(m)}^\epsilon
(\mathbf{z}_{g_{n(m)}^\epsilon}(t))\right)^2|\mathbf{z}_{g_{n(m)}^\epsilon}^\prime(t)|dt
\label{3.8}\\
&\leq &
\int_E j\left(\mathbf{x}, \widetilde{y}_{g_{m}}(\mathbf{x})\right)d\mathbf{x}
=\int_E j\left(\mathbf{x}, y_{g_{m}}(\mathbf{x})\right)d\mathbf{x}
\rightarrow \inf (\mathcal{P}).
\nonumber
\end{eqnarray}

Since $j$ is bounded from below, we get from (\ref{3.8}):
\begin{equation}\label{3.9}
\int_{\partial \Omega_{g_{m}}}\left(y_{n(m)}^\epsilon\right)^2 d\sigma \leq C\epsilon
\end{equation}
with $C$ a constant independent of $\epsilon >0$. Then, (\ref{3.9}) shows that 
(\ref{1.3}) is fulfilled with a perturbation of order $\epsilon^{1/2}$.

Moreover, again by (\ref{3.8}), we see the minimizing property of $\{y_{n(m)}^\epsilon\}$
in the original problem $(\mathcal{P})$.

We notice that in the state equation (\ref{3.3}), the right-hand side coincides with $f$
in the set
$\left\{ (x_1,x_2)\in D;\ g(x_1,x_2) < -\epsilon\right\},$
which is an approximation of $\Omega_{g_{n(m)}^\epsilon}$. Namely, we notice that for any
$g\in\mathcal{F}$, the open sets 
$\left\{ (x_1,x_2)\in D;\ g(x_1,x_2) < -\epsilon\right\}$
form a nondecreasing sequence contained in $\overline{\Omega}_g$, when $\epsilon\rightarrow 0$.
Take $(x_1,x_2)$ such that $g(x_1,x_2)=0$ and take some sequence
$(x_1^n,x_2^n)\rightarrow (x_1,x_2)$,
$(x_1^n,x_2^n)\in \Omega_g$. We have $g(x_1^n,x_2^n)<0$  by (\ref{3.1}) and $g(x_1^n,x_2^n) \rightarrow 0$. Moreover, $(x_1^n,x_2^n)\in \Omega_{\epsilon_n}=
\Omega_{g+\epsilon_n}$, for $\epsilon_n>0$ sufficiently small.
Consequently, we have the desired convergence property by \cite{Tiba2006}, p. 461.
This ends the proof.
\quad$\Box$

\begin{remark}\label{rem:3}
  A detailed study of the approximation properties in the penalized problem is performed
  in \cite{Tiba2018a}, in a slightly different case.
\end{remark}

We consider now variations $u+\lambda v$, $g+\lambda r$, where $\lambda\in \mathbb{R}$,
$u,v\in L^2(D)$,
$g,r\in\mathcal{F}$, $g(x_1^0,x_2^0)=r(x_1^0,x_2^0)=0$. Notice that  $u+\lambda v \in L^2(D)$ and
 $g+\lambda r\in \mathcal{F}$ for $|\lambda|$ sufficiently small . The conditions  (\ref{2.6}), (\ref{2.7}), (\ref{2.8}) from the
definition
of $\mathcal{F}$ are satisfied for $|\lambda|$ sufficiently small (depending on $g$) due to
the Weierstrass theorem and the fact that $\overline{E}$, $\partial D$ and $\mathcal{G}$ are compacts.
Here, we also use 
Proposition \ref{prop:2.3}. Consequently, we assume $|\lambda|$  ``small''.
We study first the differentiability properties of the state system (\ref{3.3})-(\ref{3.7}):

\begin{proposition}\label{prop:3.2}
The system of variations corresponding to (\ref{3.3})-(\ref{3.7}) is
\begin{eqnarray}
-\Delta q_\epsilon & = & (g+\epsilon)_+^2v + 2(g+\epsilon)_+u\,r,\quad\hbox{in } D,
\label{3.10}\\
q_\epsilon & = & 0,\quad\hbox{on } \partial D,
\label{3.11}\\
w_1^\prime& = & -\nabla\partial_2 g(\mathbf{z}_g)\cdot \mathbf{w}
-\partial_2 r(\mathbf{z}_g),\quad\hbox{in } I_g,
\label{3.12}\\
w_2^\prime& = & \nabla\partial_1 g(\mathbf{z}_g)\cdot \mathbf{w}
+\partial_1 r(\mathbf{z}_g),\quad\hbox{in } I_g,
\label{3.13}\\
w_1(0)& = &0,\ w_2(0) = 0,
\label{3.14}
\end{eqnarray}
where $q_\epsilon=\lim_{\lambda\rightarrow 0}\frac{y_\epsilon^\lambda-y_\epsilon}{\lambda}$, 
$\mathbf{w}=[w_1,w_2]=\lim_{\lambda\rightarrow 0} \frac{\mathbf{z}_{g+\lambda r} - \mathbf{z}_g}{\lambda}$
with $y_\epsilon^\lambda\in H^2(D)\cap H^1_0(D)$ being the solution of (\ref{3.3})-(\ref{3.4})
corresponding to $g+\lambda r$, $u+\lambda v$, and $\mathbf{z}_{g+\lambda r}\in \mathcal{C}^1(I_g)^2$
is the solution of (\ref{3.5})-(\ref{3.7}) corresponding to $g+\lambda r$. The limits exist
in the above spaces. We denote by ``$\cdot$'' the scalar product on $\mathbb{R}^2$.
\end{proposition}

\noindent
\textbf{Proof.} 
We subtract the equations corresponding to $y_\epsilon^\lambda$ and $y_\epsilon$ and divide
by $\lambda\neq 0$, small:
\begin{equation}\label{3.15}
  -\Delta  \frac{y_\epsilon^\lambda-y_\epsilon}{\lambda}
  = \frac{1}{\lambda}\left[ (g+\lambda r+\epsilon)_+^2(u+\lambda v)
  - (g+\epsilon)_+^2 u\right],\quad\hbox{in } D,
\end{equation}
with 0 boundary conditions on $\partial D$.
The regularity conditions on $\mathcal{F}$ and $u,v\in L^2(D)$ give the convergence
of the right-hand side in (\ref{3.15}) to the right-hand side in (\ref{3.10})
(strongly in $L^2(D)$) via some calculations. Then, by elliptic regularity,
we have $\frac{y_\epsilon^\lambda-y_\epsilon}{\lambda}\rightarrow q_\epsilon$ strongly in
$H^2(D)\cap H^1_0(D)$ and (\ref{3.10}), (\ref{3.11}) follows.

For (\ref{3.12})-(\ref{3.14}), the argument is the same as in Proposition 6, \cite{Tiba2013}.
The convergence of the ratio $\frac{\mathbf{z}_{g+\lambda r} - \mathbf{z}_g}{\lambda}$
is in $\mathcal{C}^1(I_g)^2$
on the whole sequence $\lambda\rightarrow 0$, due to the uniqueness property of the
linear system (\ref{3.12})-(\ref{3.14}).
Here, we also use Remark \ref{rem:2}, on the convergence $\mathcal{G}_\lambda \rightarrow \mathcal{G}$ and the
continuity with respect to the perturbations of $g$ in the Hamiltonian system
(\ref{2.1})-(\ref{2.3}), according to \cite{Tiba2013}.
\quad$\Box$

\begin{remark}\label{rem:4}
  We have as well imposed the condition
\begin{equation}\label{3.16}
g(x_1^0,x_2^0)=0,\quad \forall g \in \mathcal{F},
\end{equation}
where $(x_1^0,x_2^0)\in D\setminus E$ is some given point.
Similarly, constraints like (\ref{3.16}) may be imposed on a finite number of points
or on some curves in $D\setminus E$ and their geometric meaning is that the boundary
$\partial \Omega_g$ of the admissible unknown domains should contain these points, curves, etc.
\end{remark}

\begin{proposition}\label{prop:3.3}
Assume that $f\in L^p(D)$, $j(\mathbf{x},\cdot)$ is of class $\mathcal{C}^1(\mathbb{R})$
and bounded, $g \in \mathcal{F}$, $u\in L^p(D)$, $p>2$ and
$y_\epsilon (z_g (t)) = 0$ in $[0,T_g]$. Then, for any
direction $[r,v]\in \mathcal{F}\times L^p(D)$, the derivative of the penalized cost
(\ref{3.2}) is given by:
\begin{eqnarray}
&&
\quad \int_E \partial_2 j\left(\mathbf{x}, y_\epsilon(\mathbf{x})\right)q_\epsilon(\mathbf{x}) d\mathbf{x}
+\frac{2}{\epsilon}\int_{I_g} y_\epsilon(\mathbf{z}_g(t))
q_\epsilon(\mathbf{z}_g(t)) |\mathbf{z}_g^\prime(t)| dt
\label{3.17}\\
& + & \frac{2}{\epsilon}\int_{I_g} y_\epsilon(\mathbf{z}_g(t))\nabla y_\epsilon(\mathbf{z}_g(t))
\cdot \mathbf{w}(t) |\mathbf{z}_g^\prime(t)| dt
+\frac{1}{\epsilon}\int_{I_g} \left(y_\epsilon(\mathbf{z}_g(t))\right)^2
\frac{\mathbf{z}_g^\prime(t)\cdot \mathbf{w}^\prime(t)}{|\mathbf{z}_g^\prime(t)|} dt
\nonumber
\end{eqnarray}
where $q_\epsilon\in W^{2,p}(D)\cap W_0^{1,p}(D)$, $\mathbf{w}\in \mathcal{C}^1(I_g)^2$,
$\mathbf{z}_g\in \mathcal{C}^1(I_g)^2$ satisfy (\ref{3.10})-(\ref{3.14})
and (\ref{2.1})-(\ref{2.3})
respectively, and $I_g=[0,T_g]$ is the period interval for $\mathbf{z}_g(\cdot)$.
\end{proposition}

\noindent
\textbf{Proof.} 
In the notations of Proposition \ref{prop:3.2}, we compute
\begin{eqnarray}
  &&\lim_{\lambda \rightarrow 0}
  \left\{ \frac{1}{\lambda}
  \int_E  \left[ j\left(\mathbf{x}, y_\epsilon^\lambda(\mathbf{x})\right)
                -j\left(\mathbf{x}, y_\epsilon(\mathbf{x})\right)
          \right] d\mathbf{x}
  \right.
  \label{3.18}\\
  &&
\left.
  +\frac{1}{\epsilon\lambda}
  \int_{I_g}
  \left[
    \left(y_\epsilon^\lambda(\mathbf{z}_{g+\lambda r}(t)\right)^2 |\mathbf{z}_{g+\lambda h}^\prime(t)| 
   -\left(y_\epsilon(\mathbf{z}_{g}(t)\right)^2 |\mathbf{z}_g^\prime(t)| 
   \right] dt 
\right\}  .
\nonumber
\end{eqnarray}

In (\ref{3.18}), $\lambda >0$ is ``small'' and Proposition \ref{prop:2.3} ensures that
$g+\lambda r\in \mathcal{F}$ (see \cite{Tiba2018} as well). By Proposition \ref{prop:2.2}
we know that the trajectories associated to $g+\lambda h$ are periodic, that is the
functions in the second integrals are defined on $I_g$.
Moreover, since $f,u\in L^p(D)$, then $y_\epsilon^\lambda,\ y_\epsilon$ defined as in 
(\ref{3.3}), (\ref{3.4}) are in $W^{2,p}(D)\subset \mathcal{C}^1(\overline{D})$, by the Sobolev
theorem and elliptic regularity. Consequently, all the integrals
appearing in (\ref{3.17}), (\ref{3.18}) make sense.

Moreover, in (\ref{3.18}), we have neglected the term
\begin{eqnarray*}
  L & = & \lim_{\lambda\rightarrow 0}\frac{1}{\lambda\epsilon}
  \int_{T_g}^{T_{g+\lambda r}} y_\epsilon^\lambda
  \left(\mathbf{z}_{g+\lambda r}(t)\right)^2 |\mathbf{z}_{g+\lambda r}^\prime(t)| dt\\
  &=& \lim_{\lambda\rightarrow 0}\frac{1}{\lambda\epsilon}
  \int_{T_g}^{T_{g+\lambda r}}\left[
  y_\epsilon^\lambda\left(\mathbf{z}_{g+\lambda r}(t)\right)^2
  |\mathbf{z}_{g+\lambda r}^\prime(t)|
  -y_\epsilon\left(\mathbf{z}_{g}(t)\right)^2
  |\mathbf{z}_{g}^\prime(t)|
  \right]dt
\end{eqnarray*}
due to the hypothesis on $y_\epsilon\left(\mathbf{z}_{g}(t)\right)$. We can study term by term:
\begin{eqnarray*}
  \int_{T_g}^{T_{g+\lambda r}}
  \left[
    \frac{y_\epsilon^\lambda\left(\mathbf{z}_{g+\lambda r}(t)\right)^2
   -y_\epsilon\left(\mathbf{z}_{g+\lambda r}(t)\right)^2}{\lambda}
  |\mathbf{z}_{g+\lambda r}^\prime(t)|
  \right] dt,\\
\int_{T_g}^{T_{g+\lambda r}}
  \left[
    \frac{y_\epsilon\left(\mathbf{z}_{g+\lambda r}(t)\right)^2
   -y_\epsilon\left(\mathbf{z}_{g}(t)\right)^2}{\lambda}
  |\mathbf{z}_{g+\lambda r}^\prime(t)|
  \right] dt,\\
  \int_{T_g}^{T_{g+\lambda r}}
  \left[
    y_\epsilon\left(\mathbf{z}_{g}(t)\right)^2
    \frac{|\mathbf{z}_{g+\lambda r}^\prime(t)|
    -|\mathbf{z}_{g}^\prime(t)|}{\lambda}
  \right] dt .
\end{eqnarray*}

By Proposition \ref{prop:3.2}, each of the above three
integrands are uniformly bounded and their limits can be
easily computed, for instance on
$[0,2T_g]$ due to Proposition \ref{prop:2.4}.
Notice, in the last term, that
$|\mathbf{z}_{g}^\prime(t)|
=|\nabla g\left(\mathbf{z}_{g}(t)\right)|\neq 0$ due to
(\ref{3.5}), (\ref{3.6}) and (\ref{3.7}), that is we can
differentiate here as well.

Again by Proposition \ref{prop:2.4} and the above uniform
boundedness, we infer that each of the above three terms
has null limit as $\lambda \rightarrow 0$, i.e. $L=0$.
Consequently, it is enough to study the limit (\ref{3.18}).

We also have $y_\epsilon^\lambda \rightarrow y_\epsilon$ in $\mathcal{C}^1(\overline{D})$, for
$\lambda\rightarrow 0$, by elliptic regularity. Then, under the assumptions on $j(\cdot,\cdot)$,
we get
\begin{equation}\label{3.19}
\frac{1}{\lambda}
  \int_E  \left[ j\left(\mathbf{x}, y_\epsilon^\lambda(\mathbf{x})\right)
                -j\left(\mathbf{x}, y_\epsilon(\mathbf{x})\right)
          \right] d\mathbf{x}
\rightarrow
\int_E \partial_2 j\left(\mathbf{x}, y_\epsilon(\mathbf{x})\right)q_\epsilon(\mathbf{x}) d\mathbf{x} .
\end{equation}

For the second integral in (\ref{3.18}), we intercalate certain intermediary terms and we
compute their limits for $\lambda \rightarrow 0$:
\begin{eqnarray}
&&
\lim_{\lambda \rightarrow 0}
\frac{1}{\epsilon\lambda}
  \int_{I_g}
  \left[
    \left(y_\epsilon^\lambda(\mathbf{z}_{g+\lambda r}(t)\right)^2 |\mathbf{z}_{g+\lambda r}^\prime(t)| 
   -\left(y_\epsilon(\mathbf{z}_{g+\lambda r}(t)\right)^2 |\mathbf{z}_{g+\lambda r}^\prime(t)| 
   \right] dt 
\label{3.20}
\\
&=& \frac{2}{\epsilon}\int_{I_g} y_\epsilon(\mathbf{z}_g(t)) q_\epsilon(\mathbf{z}_g(t))
 |\mathbf{z}_g^\prime(t)| dt 
\nonumber
\end{eqnarray}
due to the convergence $\mathbf{z}_{g+\lambda r}\rightarrow\mathbf{z}_g$ in  $\mathcal{C}^1(I_g)^2$
by $g,r\in \mathcal{C}^2(\overline{D})$ and the continuity properties in (\ref{2.1})-(\ref{2.3});
\begin{eqnarray}
&&
\lim_{\lambda \rightarrow 0}
\frac{1}{\epsilon\lambda}
  \int_{I_g}
  \left[
    \left(y_\epsilon(\mathbf{z}_{g+\lambda r}(t)\right)^2 |\mathbf{z}_{g+\lambda r}^\prime(t)| 
   -\left(y_\epsilon(\mathbf{z}_{g}(t)\right)^2 |\mathbf{z}_{g+\lambda r}^\prime(t)| 
   \right] dt 
\label{3.21}
\\
&=& \frac{2}{\epsilon}\int_{I_g} y_\epsilon(\mathbf{z}_g(t))\nabla y_\epsilon(\mathbf{z}_g(t))
\cdot \mathbf{w}(t) |\mathbf{z}_g^\prime(t)| dt ,
\nonumber
\end{eqnarray}
where $\mathbf{w}=(w_1,w_2)$ satisfies (\ref{3.12})-(\ref{3.14}) and again we use the regularity and the
convergence properties in $\mathcal{C}^1(D)$, respectively $\mathcal{C}^1(I_g)^2$.
\begin{eqnarray}
&&
\lim_{\lambda \rightarrow 0}
\frac{1}{\epsilon\lambda}
  \int_{I_g}
  \left[
    \left(y_\epsilon(\mathbf{z}_{g}(t)\right)^2 |\mathbf{z}_{g+\lambda r}^\prime(t)| 
   -\left(y_\epsilon(\mathbf{z}_{g}(t)\right)^2 |\mathbf{z}_{g}^\prime(t)| 
   \right] dt 
\label{3.22}
\\
&=& \frac{1}{\epsilon}
\int_{I_g} \left(y_\epsilon(\mathbf{z}_g(t))\right)^2
\frac{\mathbf{z}_g^\prime(t)\cdot \mathbf{w}^\prime(t)}{|\mathbf{z}_g^\prime(t)|} dt,
\nonumber
\end{eqnarray}
where we recall that $|\mathbf{z}_g^\prime(t)|=\sqrt{(z_g^1)^\prime(t)^2+(z_g^2)^\prime(t)^2}$ is non zero by (\ref{2.7}) 
and the Hamiltonian system,
and standard derivation rules may be applied under our regularity conditions.

By summing up (\ref{3.19})-(\ref{3.22}), we end the proof of (\ref{3.17}).
\quad$\Box$

\begin{remark}\label{rem:5}
In the case that $\Omega_g$ is not simply connected, the penalization integral in
(\ref{3.2}) is in fact a finite sum and each of these terms can be handled separately,
in the same way as above, due to Proposition \ref{prop:2.3} and Remark \ref{rem:2}.
The significance of the hypothesis
$y_\epsilon(\mathbf{z}_g(t))=0$, is that one should first
minimize the penalization term with respect to the control
$u$ (this is possible due to the arguments in the proof
of Proposition \ref{prop:3.1}). Then, the obtained control
should be fixed and the minimization with respect to
$g\in\mathcal{F}$ is to be performed.
In case $T_{g+\lambda r}$ can be evaluated as in
Remark \ref{rem:2.7}, then the hypothesis can be relaxed 
to $y_\epsilon(x_1^0,x_2^0)=0$ via a variant of the above
arguments.

\end{remark}

Now, we denote by
$A:\mathcal{C}^2(\overline{D})\times L^p(D)\rightarrow W^{2,p}(D)\cap W_0^{1,p}(D)$
the linear continuous operator given by
$r, v\rightarrow q_\epsilon$, defined in (\ref{3.10}), (\ref{3.11}).
We also denote by $B: \mathcal{C}^2(\overline{D}) \rightarrow \mathcal{C}^1(I_g)^2$
the linear continuous operator given by (\ref{3.12})-(\ref{3.14}),
$Br=[w_1,w_2]$. In these definitions, $g\in \mathcal{C}^2(\overline{D})$ and
$u\in L^p(D)$ are fixed. We have:

\begin{corollary}\label{cor:1}
The relation (\ref{3.17}) can be rewritten as:
\begin{eqnarray}
&&
\int_E \partial_2 j\left(\mathbf{x}, y_\epsilon(\mathbf{x})\right)A(r,v)(\mathbf{x}) d\mathbf{x}
+\frac{2}{\epsilon}\int_{I_g} y_\epsilon(\mathbf{z}_g(t))
A(r,v)(\mathbf{z}_g(t)) |\mathbf{z}_g^\prime(t)| dt
\label{3.23}\\
&& +\frac{2}{\epsilon}\int_{I_g} y_\epsilon(\mathbf{z}_g(t))\nabla y_\epsilon(\mathbf{z}_g(t))
\cdot Br(t) |\mathbf{z}_g^\prime(t)| dt
\nonumber\\
&&
+\frac{1}{\epsilon}\int_{I_g}
\frac{\left(y_\epsilon(\mathbf{z}_g(t))\right)^2}{|\mathbf{z}_g^\prime(t)|}
\mathbf{z}_g^\prime(t)\cdot [-\partial_2 r,\partial_1 r](\mathbf{z}_g(t)) dt
\nonumber\\
&&
+\frac{1}{\epsilon}\int_{I_g}
\frac{\left(y_\epsilon(\mathbf{z}_g(t))\right)^2}{|\mathbf{z}_g^\prime(t)|}
C(t) \cdot \mathbf{w}(t) dt ,
\nonumber
\end{eqnarray}
where the vector $C(t)$ is explained below.
\end{corollary}

\noindent
\textbf{Proof.}
In the last integral in (\ref{3.17}), we replace $\mathbf{w}^\prime(t)$ by the
right-hand side in (\ref{3.12}), (\ref{3.13}). We compute:
\begin{eqnarray}
&&  \mathbf{z}_g^\prime(t) \cdot \mathbf{w}^\prime(t)
\label{3.24}
\\
&=&
\mathbf{z}_g^\prime(t) \cdot
     [-\nabla\partial_2 g(\mathbf{z}_g(t))\cdot \mathbf{w}(t)-\partial_2 r(\mathbf{z}_g(t)),
      \nabla\partial_1 g(\mathbf{z}_g(t))\cdot \mathbf{w}(t) +\partial_1 r(\mathbf{z}_g(t))
     ]  
\nonumber
\\
&=&
\mathbf{z}_g^\prime(t) \cdot[-\partial_2 r(\mathbf{z}_g(t)),\partial_1 r(\mathbf{z}_g(t))]  
\nonumber
\\&+&
\mathbf{z}_g^\prime(t) \cdot
     [-\partial_{1,2}^2 g(\mathbf{z}_g(t)) w_1(t),
      \partial_{1,1}^2 g(\mathbf{z}_g(t)) w_1(t)
     ]
     \nonumber
\\&+&
\mathbf{z}_g^\prime(t) \cdot
     [
      -\partial_{2,2}^2 g(\mathbf{z}_g(t)) w_2(t),
      \partial_{2,1}^2 g(\mathbf{z}_g(t)) w_2(t) 
     ].
\nonumber
\end{eqnarray}
We denote by $C(t)$ the (known) vector
\begin{eqnarray*}
C(t)&=&[-(z_g^1)^\prime(t)\partial_{1,2}^2 g(\mathbf{z}_g(t))
        +(z_g^2)^\prime(t)\partial_{1,1}^2 g(\mathbf{z}_g(t))  ,
\\
&&
-(z_g^1)^\prime(t)\partial_{2,2}^2 g(\mathbf{z}_g(t)) 
+(z_g^2)^\prime(t)\partial_{2,1}^2 g(\mathbf{z}_g(t)) 
     ]
\end{eqnarray*}
and together with (\ref{3.24}), we get (\ref{3.23}). This ends the proof.
\quad$\Box$

\section{Finite element discretization}
\setcounter{equation}{0}

We assume that $D$ and $E$ are polygonal. Let $\mathcal{T}_h$ be a triangulation of $D$ with vertices $A_i$,
$i\in I=\{1,\dots,n\}$. We consider that $\mathcal{T}_h$ is compatible with $E$, i.e. 
$$
\forall T\in \mathcal{T}_h, T \subset \overline{E} \hbox{ or } T\subset \overline{D\setminus E}
$$
where $T$ designs a triangle of $\mathcal{T}_h$ and $h$ is the size of $\mathcal{T}_h$. 
We consider a triangle as a closed set.
For simplicity, we employ piecewise linear finite element and we denote
$$
\mathbb{W}_h=\{ \varphi_h\in \mathcal{C}(\overline{D});
\ {\varphi_h}_{|T} \in \mathbb{P}_1(T),\ \forall T \in \mathcal{T}_h  \}.
$$
We use a standard basis of $\mathbb{W}_h$, $\{ \phi_i\}_{i\in I}$, where $\phi_i$ is the hat function
associated to the vertex $A_i$, see for example \cite{Ciarlet2002}, \cite{Raviart2004}. 
The finite element approximations of $g$ and $u$ are
$g_h(\mathbf{x})=\sum_{i\in I} G_i \phi_i(\mathbf{x})$, 
$u_h(\mathbf{x})=\sum_{i\in I} U_i \phi_i(\mathbf{x})$,
for all $\mathbf{x}\in \overline{D}$.  We set the vectors
$G=(G_i)_{i\in I}^T\in\mathbb{R}^n$, $U=(U_i)_{i\in I}^T\in\mathbb{R}^n$
and $g_h$ can be identified by $G$, etc.
The function $u$ is in $L^p(D)$, as in Proposition \ref{prop:3.3}.
Alternatively, for $u_h$, we can use discontinuous piecewise constant finite element $\mathbb{P}_0$.
In order to approach $g\in \mathcal{C}^2(\overline{D})$, we can use high order finite elements.

\subsection{Discretization of the optimization problem}

We introduce
$$
\mathbb{V}_h=\{ \varphi_h\in \mathbb{W}_h;
\ \varphi_h = 0\hbox{ on }\partial D \},
$$
$I_0=\{i\in I;\ A_i\notin \partial D \}$ and $n_0=card(I_0)$.
The finite element weak formulation of (\ref{3.3})-(\ref{3.4}) is: 
find $y_h\in \mathbb{V}_h$ such that
\begin{equation}\label{4.1}
  \int_D \nabla y_h \cdot \nabla  \varphi_h d\mathbf{x}
= \int_D \left(f + (g_h+\epsilon)_+^2u_h\right) \varphi_h d\mathbf{x},\quad \forall \varphi_h \in \mathbb{V}_h.
\end{equation}
As before, for $y_h(\mathbf{x})=\sum_{j\in I_0} Y_j \phi_j(\mathbf{x})$, 
we set $Y=(Y_j)_{j\in I_0}^T\in\mathbb{R}^{n_0}$.
In order to obtain the linear system, we take the basis functions $\varphi_h=\phi_i$ in (\ref{4.1}) for 
$i\in I_0$.
Let us consider the vector
$$
F=\left( \int_D f\phi_i d\mathbf{x}\right)_{i\in I_0}^T\in\mathbb{R}^{n_0},
$$
the $n_0\times n_0$ matrix $K$ defined by
$$
K=(K_{ij})_{i\in I_0,j\in I_0},\quad K_{ij}=  \int_D \nabla  \phi_j \cdot \nabla  \phi_i d\mathbf{x}
$$
and the $n_0\times n$ matrix $B^1(G,\epsilon)$ defined by
$$
B^1(G,\epsilon)=(B_{ij}^1)_{i\in I_0,j\in I},\quad 
B_{ij}^1=  \int_D  (g_h+\epsilon)_+^2 \phi_j\phi_i d\mathbf{x} .
$$
The matrix $K$ is symmetric, positive definite.
The finite element approximation of the state system (\ref{3.3})-(\ref{3.4}) is the linear system:
\begin{equation}\label{4.2}
KY=F+B^1(G,\epsilon)U.
\end{equation}

Now, we shall discretize the objective function (\ref{3.2}).
We denote $I_E=\{i\in I;\ A_i\in \overline{E} \}$ and $n_E=card(I_E)$.
For the first term of (\ref{3.2}), we introduce
$$
J_1(Y)=\int_E j(\mathbf{x}, y_h(\mathbf{x})) d\mathbf{x}.
$$

We shall study the second term of (\ref{3.2}).
In order to solve numerically the ODE system (\ref{3.5})-(\ref{3.7}), we use a
partition $[t_0,\dots,t_k,\dots,t_m]$ of $[0,T_g]$, with $t_0=0$ and $t_m=T_g$.
We can use the forward Euler scheme:
\begin{eqnarray}
Z_{k+1}^1 & = & Z_{k}^1 - (t_{k+1}-t_k) \frac{\partial g_h}{\partial x_2}\left(Z_{k}^1,Z_{k}^2\right),
\label{4.3}\\
Z_{k+1}^2 & = & Z_{k}^2 + (t_{k+1}-t_k) \frac{\partial g_h}{\partial x_1}\left(Z_{k}^1,Z_{k}^2\right),
\label{4.4}\\
(Z_0^1,Z_0^2) & = & \left(x_1^0,x_2^0\right),
\label{4.5}
\end{eqnarray}
for $k=0,\dots , m-2$. We set $Z_k=(Z_k^1,Z_k^2)$ and we impose $Z_m=Z_0$.
In fact, $Z_k$ is an approximation of $\mathbf{z}_g(t_k)$.
We do not need to stock $Z_0$ and we set $Z=(Z^1,Z^2)\in \mathbb{R}^m\times \mathbb{R}^m$, with
$Z^1=(Z_k^1)_{1\leq k\leq m}^T$ and $Z^2=(Z_k^2)_{1\leq k\leq m}^T$.
In the applications, one can use more performant numerical methods for the ODE's,
like explicit Runge-Kutta or backward Euler, 
but here we want to avoid a too tedious exposition.

Without risk of confusion, we introduce the function $Z:[0,T_g]\rightarrow \mathbb{R}^2$ defined by
$$
Z(t)=\frac{t_{k+1}-t}{(t_{k+1}-t_k)} Z_k + \frac{t-t_k}{(t_{k+1}-t_k)} Z_{k+1},\quad t_k \leq t < t_{k+1}
$$
for $k=0,1,\dots ,m-1$. We have $Z(t_k)=Z_k$ and we can identify the function $Z(\cdot)$ by the vector 
$Z\in \mathbb{R}^m\times \mathbb{R}^m$. We remark that $Z(\cdot)$ is derivable on each interval $(t_k, t_{k+1})$
and $Z^\prime(t)=\frac{1}{(t_{k+1}-t_k)}(Z_{k+1}^1-Z_k^1,Z_{k+1}^2-Z_k^2)$ for $t_k \leq t < t_{k+1}$.

We introduce the $n_0\times n_0$ matrix $N(Z)$ defined by
$$
N(Z)=\left( 
\int_0^{T_g} \phi_j(Z(t))  \phi_i(Z(t)) |Z^\prime(t)| dt
\right)_{i\in I_0,j\in I_0}
$$
and the second term of (\ref{3.2}) is approached by
$\frac{1}{\epsilon}Y^TN(Z)Y$, then the discrete form of the optimization problem (\ref{3.2})-(\ref{3.7})
is
\begin{equation}\label{4.6}
\min_{G,U\in \mathbb{R}^n} J(G,U)=J_1(Y)+\frac{1}{\epsilon}Y^TN(Z)Y
\end{equation}
subject to (\ref{4.2}). We point out that $Y$ depends on $G$ and $U$ by (\ref{4.2})
and $Z$ depends on $G$ by (\ref{4.3})-(\ref{4.5}).

\subsection{Discretization of the derivative of the objective function}

Let $r_h,v_h$ be in $\mathbb{W}_h$ and $R,V\in \mathbb{R}^n$ be the associated vectors.
The finite element weak formulation of (\ref{3.10})-(\ref{3.11}) is: 
find $q_h\in \mathbb{V}_h$ such that
\begin{equation}\label{4.7}
  \int_D \nabla q_h \cdot \nabla  \varphi_h d\mathbf{x}
  = \int_D \left( (g_h+\epsilon)_+^2 v_h + 2(g_h+\epsilon)_+ u_h r_h\right)
  \varphi_h d\mathbf{x},
  \quad \forall \varphi_h \in \mathbb{V}_h.
\end{equation}
Let $Q\in \mathbb{R}^{n_0}$ be the associated vector to $q_h$ and we construct
the $n_0\times n$ matrix $C^1(G,\epsilon,U)$ defined by
$$
C^1(G,\epsilon,U)=\left(
\int_D  2(g_h+\epsilon)_+ u_h\phi_j\phi_i d\mathbf{x}
\right)_{i\in I_0,j\in I}.
$$
The linear system of (\ref{4.7}) is
\begin{equation}\label{4.8}
KQ=B^1(G,\epsilon)V + C^1(G,\epsilon,U)R.
\end{equation}

In order to approximate $\partial_2 j(\mathbf{x}, y_\epsilon(\mathbf{x}))$, $y_\epsilon$ 
given by (\ref{3.3})-(\ref{3.4}),
we consider the nonlinear application
$$
Y\in \mathbb{R}^{n_0} \rightarrow L(Y)\in \mathbb{R}^{n_E}
$$
such that
$\partial_2 j(\mathbf{x}, y_h(\mathbf{x}))
=\sum_{i\in I_E} \left(L(Y)\right)_i {\phi_i}_{|E} (\mathbf{x})$
where ${\phi_i}_{|E}$ is the restriction of $\phi_i$ to $E$.
We define the $n_E\times n_0$ matrix $M_{ED}$ defined by
$$
M_{ED}=\left(
\int_D \phi_i\phi_j d\mathbf{x}
\right)_{i\in I_E,j\in I_0}.
$$
The first term of (\ref{3.17}) is approached by
\begin{equation}\label{4.9}
  \left(L(Y)\right)^T M_{ED} Q
\end{equation}
and the second term of (\ref{3.17}) is approached by
\begin{equation}\label{4.10}
\frac{2}{\epsilon}Y^T N(Z)Q
\end{equation}
where the matrix $N(Z)$ was introduced in the previous subsection.

Next, we introduce the partial derivative for a piecewise linear function.
Let $g_h\in \mathbb{W}_h$ and $G\in \mathbb{R}^n$ its associated vector, i.e.
$g_h(\mathbf{x})=\sum_{i\in I} G_i \phi_i(\mathbf{x})$.
Let $\Pi_h^1G \in \mathbb{R}^n$ defined by
$$
\left(\Pi_h^1 G\right)_i
=\frac{1}{\sum_{j\in J_i} area(T_j)}
\sum_{j\in J_i} area(T_j)\partial_1 {g_h}_{|T_j}
$$
where $J_i$ is the set of index $j$ such that the triangle $T_j$
has the vertex $A_i$.
Since $g_h$ is a linear function in each triangle $T_j$, then $\partial_1 {g_h}_{|T_j}$
is constant.
Similarly, we construct $\Pi_h^2 G \in \mathbb{R}^n$ for $\partial_2$.
In fact, $\Pi_h^1$ and $\Pi_h^2$ are two $n\times n$ matrices depending on $\mathcal{T}_h$.

Then, we set
$$
\partial_1^h g_h(\mathbf{x})=\sum_{i\in I} \left(\Pi_h^1G\right)_i \phi_i(\mathbf{x})
$$
and similarly for $\partial_2^h g_h$.
Finally, we put $\nabla_h g_h = (\partial_1^h g_h, \partial_2^h g_h)$.
Since $y_h\in \mathbb{V}_h\subset \mathbb{W}_h$, we can define $\partial_1^h y_h$
and $\partial_2^h y_h$.

\begin{example}
We shall give a simple example to understand the discrete derivative of $\mathbb{W}_h$ functions.
We consider the square $[A_1A_2A_4A_3]$ of vertices $A_1=(0,0)$, $A_2=(1,0)$, $A_4=(1,1)$,
$A_3=(0,1)$ and the triangulation of two triangles $T_1=[A_1A_2A_4]$ and $T_2=[A_1A_4A_3]$.
We shall present the discrete derivative of the hat function
$$
\phi_4(x_1,x_2) =
\left\{
\begin{array}{ll}
  x_2 &\hbox{ in }T_1\\
  x_1 &\hbox{ in }T_2 .
\end{array}  
\right.
$$
We have $J_1=\{1,2\}$ and 
\begin{eqnarray*}
\left(\Pi_h^1 \phi_4 \right)_1
&=&\frac{1}{area(T_1)+area(T_2)}
\left( area(T_1)\partial_1 {\phi_4}_{|T_1}
+area(T_2)\partial_1 {\phi_4}_{|T_2}
\right)\\
&=&\frac{1}{1/2+1/2}
\left( 1/2 \times 0
+1/2 \times 1
\right)
=1/2.
\end{eqnarray*}
Similarly, $J_2=\{1\}$, $J_3=\{2\}$, $J_4=\{1,2\}$,
\begin{eqnarray*}
\left(\Pi_h^1 \phi_4 \right)_2
&=&\frac{1}{1/2}
\left( 1/2 \times 0
\right)
=0,\\
\left(\Pi_h^1 \phi_4 \right)_3
&=&\frac{1}{1/2}
\left( 1/2 \times 1
\right)
=1/2,\\
\left(\Pi_h^1 \phi_4 \right)_4
&=&\frac{1}{1/2+1/2}
\left( 1/2 \times 0
+1/2 \times 1
\right)
=1/2
\end{eqnarray*}
then
$$
\partial_1^h \phi_4(x_1,x_2)= 1/2 \times \phi_1(x_1,x_2)+ 0\times \phi_2(x_1,x_2)
+ 1 \times \phi_3(x_1,x_2)+ 1/2 \times \phi_4(x_1,x_2).
$$
\end{example}

In order to solve the ODE system (\ref{3.12})-(\ref{3.14}), we use
the forward Euler scheme on the same partition as for (\ref{4.3})-(\ref{4.5}):
\begin{eqnarray}
W_{k+1}^1  & = & W_{k}^1
  - (t_{k+1}-t_k) \nabla_h \partial_2^h g_h(Z_{k})\cdot (W_k^1,W_k^2)\label{4.11}\\
&&
- (t_{k+1}-t_k) \partial_2^h r_h(Z_{k}),
\nonumber\\
 W_{k+1}^2  &= & W_{k}^2
+ (t_{k+1}-t_k) \nabla_h \partial_1^h g_h \left(Z_{k}\right)\cdot (W_k^1,W_k^2)\label{4.12}\\
&& + (t_{k+1}-t_k) \partial_1^h r_h\left(Z_{k}\right),
\nonumber\\
 W_0^1 & = & 0,\ W_0^2=0,
\label{4.13}
\end{eqnarray}
for $k=0,\dots , m-1$.
We set $W_k=(W_k^1,W_k^2)$ and now we have $W_m \neq W_0$ generally.
In fact, $W_k$ is an approximation of $\mathbf{w}(t_k)$.
We do not need to stock $W_0$ and we set
$W=(W^1,W^2)\in \mathbb{R}^m\times \mathbb{R}^m$, with
$W^1=(W_k^1)_{1\leq k\leq m}^T$ and $W^2=(W_k^2)_{1\leq k\leq m}^T$.
As mentioned before, we can use more performant numerical methods for the ODE,
like explicit Runge-Kutta or backward Euler.

We construct $W:[0,T_g]\rightarrow \mathbb{R}^2$ in the same way as for $Z(t)$
$$
W(t)=\frac{t_{k+1}-t}{(t_{k+1}-t_k)} W_k + \frac{t-t_k}{(t_{k+1}-t_k)} W_{k+1},\quad t_k \leq t < t_{k+1}
$$
for $k=0,1,\dots ,m-1$. We have $W(t_k)=W_k$ and
$W^\prime(t)=\frac{1}{(t_{k+1}-t_k)}(W_{k+1}^1-W_k^1,W_{k+1}^2-W_k^2)$ for $t_k \leq t < t_{k+1}$.
If $\psi_k$ is the one-dimensional piecewise linear hat function associated to the point $t_k$ of the
partition $[t_0,\dots,t_k,\dots,t_m]$, we can write equivalently
$W(t)=\sum_{k=0}^m W_k \psi_k(t)$ for $t\in [0,T_g]$.
The third term of (\ref{3.17}) is approached by
\begin{eqnarray}
&&\frac{2}{\epsilon}\sum_{k=0}^{m-1}\int_{t_k}^{t_{k+1}} y_h(Z(t))\nabla_h y_h(Z(t))
\cdot W(t) |Z^\prime(t)| dt
\label{4.14}\\
&=& \frac{2}{\epsilon}\sum_{k=0}^{m-1}\int_{t_k}^{t_{k+1}} y_h(Z(t))\nabla_h y_h(Z(t))
\cdot \left( W_k\psi_k(t)+W_{k+1}\psi_{k+1}(t)\right) |Z^\prime(t)| dt
\nonumber\\
&=& \frac{2}{\epsilon}\sum_{k=0}^{m-1}\int_{t_k}^{t_{k+1}} y_h(Z(t))\partial_1^h y_h(Z(t))
\left( W_k^1\psi_k(t)+W_{k+1}^1\psi_{k+1}(t)\right)\frac{|Z_kZ_{k+1}|}{(t_{k+1}-t_k)} dt
\nonumber\\
&+& \frac{2}{\epsilon}\sum_{k=0}^{m-1}\int_{t_k}^{t_{k+1}} y_h(Z(t))\partial_2^h y_h(Z(t))
\left( W_k^2\psi_k(t)+W_{k+1}^2\psi_{k+1}(t)\right)\frac{|Z_kZ_{k+1}|}{(t_{k+1}-t_k)} dt
\nonumber
\end{eqnarray}
where $|Z_kZ_{k+1}|$ is the length of the segment in $\mathbb{R}^2$ with ends $Z_k$ and $Z_{k+1}$.

We have
\begin{eqnarray}
&&
\quad\int_{t_k}^{t_{k+1}}
y_h(Z(t))\partial_1^h y_h(Z(t))
\left( W_k^1\psi_k(t)+W_{k+1}^1\psi_{k+1}(t)\right)
\frac{|Z_kZ_{k+1}|}{(t_{k+1}-t_k)} dt
\label{4.15}\\
&=& W_k^1\int_{t_k}^{t_{k+1}}
\left(\sum_{i\in I_0} Y_i \phi_i(Z(t))\right)
\left(\sum_{j\in I} (\Pi_h^1Y)_j \phi_j(Z(t)) \right)
\psi_k(t)
\frac{|Z_kZ_{k+1}|}{(t_{k+1}-t_k)} dt
\nonumber\\
&+& W_{k+1}^1\int_{t_k}^{t_{k+1}}
\left(\sum_{i\in I_0} Y_i \phi_i(Z(t))\right)
\left(\sum_{j\in I} (\Pi_h^1Y)_j \phi_j(Z(t)) \right)
\psi_{k+1}(t)
\frac{|Z_kZ_{k+1}|}{(t_{k+1}-t_k)} dt.
\nonumber
\end{eqnarray}
We introduce the $n_0\times n$ matrices $N_k^{[k,k+1]}(Z)$ and $N_{k+1}^{[k,k+1]}(Z)$ defined by
\begin{eqnarray*}
N_k^{[k,k+1]}(Z)&=&\left( 
\int_{t_k}^{t_{k+1}} \phi_i(Z(t))  \phi_j(Z(t)) \psi_k(t) \frac{|Z_kZ_{k+1}|}{(t_{k+1}-t_k)} dt
\right)_{i\in I_0,j\in I}\\
N_{k+1}^{[k,k+1]}(Z)&=&\left( 
\int_{t_k}^{t_{k+1}} \phi_i(Z(t))  \phi_j(Z(t)) \psi_{k+1}(t) \frac{|Z_kZ_{k+1}|}{(t_{k+1}-t_k)} dt
\right)_{i\in I_0,j\in I}
\end{eqnarray*}
then (\ref{4.15}) can be rewritten as
$$
Y^T\left( W_k^1N_k^{[k,k+1]}(Z)+ W_{k+1}^1N_{k+1}^{[k,k+1]}(Z)\right)(\Pi_h^1Y)
$$
and finally the third term of (\ref{3.17}) is approached by
\begin{eqnarray}\label{4.16}
&&  \frac{2}{\epsilon}
  Y^T
  \sum_{k=0}^{m-1}\left( W_k^1 N_k^{[k,k+1]}(Z)+ W_{k+1}^1 N_{k+1}^{[k,k+1]}(Z)\right)(\Pi_h^1 Y)\\
&+&
  \frac{2}{\epsilon}
  Y^T
  \sum_{k=0}^{m-1}\left( W_k^2 N_k^{[k,k+1]}(Z)+ W_{k+1}^2 N_{k+1}^{[k,k+1]}(Z)\right)(\Pi_h^2 Y).
\nonumber  
\end{eqnarray}
We can introduce the linear operators $T^1(Z)$ and $T^2(Z)$ by 
\begin{eqnarray}\label{4.17}
\qquad W^1 \in \mathbb{R}^m \rightarrow T^1(Z)W^1&=&  
  \sum_{k=0}^{m-1}\left( W_k^1 N_k^{[k,k+1]}(Z)+ W_{k+1}^1 N_{k+1}^{[k,k+1]}(Z)\right)\\
\qquad W^2 \in \mathbb{R}^m \rightarrow T^2(Z)W^2&=&
  \sum_{k=0}^{m-1}\left( W_k^2 N_k^{[k,k+1]}(Z)+ W_{k+1}^2 N_{k+1}^{[k,k+1]}(Z)\right)
\nonumber  
\end{eqnarray}
then (\ref{4.16}) can be rewritten as
\begin{eqnarray}\label{4.18}
\frac{2}{\epsilon}
  Y^T \left( T^1(Z)W^1\right) (\Pi_h^1 Y)
+\frac{2}{\epsilon}
  Y^T \left(T^2(Z)W^2\right) (\Pi_h^2 Y).
\end{eqnarray}

The fourth term of term of (\ref{3.17}) is approached by
\begin{eqnarray}\label{4.19}
\qquad
\frac{1}{\epsilon}\sum_{k=0}^{m-1}\int_{t_k}^{t_{k+1}}
\left(\sum_{i\in I_0} Y_i \phi_i(Z(t))\right)
\left(\sum_{j\in I_0} Y_j \phi_j(Z(t))\right)
\frac{Z^\prime(t)\cdot W^\prime(t)}{|Z^\prime(t)|} dt.
\end{eqnarray}
But $Z^\prime(t)$ and $W^\prime(t)$ are constants for $t_k \leq t < t_{k+1}$, then
$$
\frac{Z^\prime(t)\cdot W^\prime(t)}{|Z^\prime(t)|} 
=\frac{(Z_{k+1}^1-Z_k^1,Z_{k+1}^2-Z_k^2)\cdot(W_{k+1}^1-W_k^1,W_{k+1}^2-W_k^2)}
{(t_{k+1}-t_k) |Z_kZ_{k+1}|}  
$$
where $|Z_kZ_{k+1}|$ is the length of the segment in $\mathbb{R}^2$ with ends $Z_k$ and $Z_{k+1}$.
We introduce the $n_0\times n_0$ matrix $R_k(Z)$ defined by
$$
R_k(Z)=\left( 
\int_{t_k}^{t_{k+1}} \phi_i(Z(t))  \phi_j(Z(t)) |Z^\prime(t)| dt
\right)_{i\in I_0,j\in I_0}
$$
and the linear operators $T^3(Z)$  
\begin{eqnarray}\label{4.20}
&&  \qquad W \in \mathbb{R}^m \times \mathbb{R}^m \rightarrow T^3(Z)W\\
&& T^3(Z)W =
  \sum_{k=0}^{m-1}\frac{(Z_{k+1}^1-Z_k^1,Z_{k+1}^2-Z_k^2)\cdot(W_{k+1}^1-W_k^1,W_{k+1}^2-W_k^2)}
      {|Z_kZ_{k+1}|^2} R_k(Z) .
  \nonumber
\end{eqnarray}
The (\ref{4.19}) can be rewritten as
\begin{eqnarray}\label{4.21}
 \frac{1}{\epsilon} Y^T \left(T^3(Z)W\right) Y.
\end{eqnarray}

The study of this subsection can be resumed as following:

\begin{proposition}\label{prop:dJ}
The discret version of (\ref{3.17}) is 
\begin{eqnarray}\label{4.22}
\ \qquad  dJ_{(G,U)}(R,V) & = &  \left(L(Y)\right)^T M_{ED} Q
  +\frac{2}{\epsilon}Y^T N(Z)Q\\
  &+& \frac{2}{\epsilon} Y^T \left( T^1(Z)W^1\right) (\Pi_h^1 Y)
  +\frac{2}{\epsilon} Y^T \left(T^2(Z)W^2\right) (\Pi_h^2 Y)\nonumber\\
  &+&\frac{1}{\epsilon} Y^T \left(T^3(Z)W\right) Y\nonumber
\end{eqnarray}
which represents the derivative of $J$ at $(G,U)$ in the direction $(R,V)$.
\end{proposition}

\noindent
\textbf{Proof.}
We get (\ref{4.22}) just by assembling (\ref{4.9}), (\ref{4.10}), (\ref{4.18}) and (\ref{4.21}).
\quad$\Box$

\subsection{Discretization of the formula (\ref{3.23})}

From (\ref{4.8}), we get
$$
Q=K^{-1}B^1(G,\epsilon)V + K^{-1}C^1(G,\epsilon,U)R
$$
and the discrete version of the operator $A$ in the Corollary \ref{cor:1} is
$$
(R,V)\in \mathbb{R}^n \times \mathbb{R}^n \rightarrow
A^1(R,V)=K^{-1}B^1(G,\epsilon)V + K^{-1}C^1(G,\epsilon,U)R.
$$
Replacing $Q$ in the first two terms of (\ref{4.22}), we get
\begin{eqnarray}\label{4.23}
  && \left( \left(L(Y)\right)^T M_{ED} +\frac{2}{\epsilon}Y^T N(Z)\right)
  K^{-1}B^1(G,\epsilon)V
  \\
  &+& \left( \left(L(Y)\right)^T M_{ED} +\frac{2}{\epsilon}Y^T N(Z)\right)
  K^{-1}C^1(G,\epsilon,U)R.
  \nonumber
\end{eqnarray}

We denote
\begin{eqnarray*}
  \Lambda_1(t)&=& y_h(Z(t))\nabla y_h(Z(t))|Z^\prime(t)|\\
  \Lambda_2(t)&=&\frac{\left(y_h(Z(t))\right)^2}{|Z^\prime(t)|} Z^\prime(t)\\
  \Lambda_3(t)&=&\frac{\left(y_h(Z(t))\right)^2}{|Z^\prime(t)|} C(t).
\end{eqnarray*}

The third term of (\ref{3.23}) is approached by
$$
\frac{2}{\epsilon}\int_0^{T_g} \Lambda_1(t) \cdot W(t)\, dt
$$
and using the trapezoidal quadrature formula on each sub-interval $[t_k,t_{k+1}]$, we get
\begin{equation}\label{4.24}
  \frac{1}{\epsilon}\sum_{k=0}^{m-1}(t_{k+1}-t_k)\left[
  \Lambda_1(t_k)\cdot (W_k^1,W_k^2)
  +\Lambda_1(t_{k+1})\cdot (W_{k+1}^1,W_{k+1}^2)
  \right].
\end{equation}
Similarly, for the 4th and 5th terms of (\ref{3.23}), we get
\begin{equation}\label{4.25}
  \frac{1}{2\epsilon}\sum_{k=0}^{m-1}(t_{k+1}-t_k)\left[
  \Lambda_2(t_k)\cdot [-\partial_2^h r_h,\partial_1^h r_h](Z_k)
  +\Lambda_2(t_{k+1})\cdot [-\partial_2^h r_h,\partial_1^h r_h](Z_{k+1})
  \right]
\end{equation}
and
\begin{equation}\label{4.26}
  \frac{1}{2\epsilon}\sum_{k=0}^{m-1}(t_{k+1}-t_k)\left[
  \Lambda_3(t_k)\cdot (W_k^1,W_k^2)
  +\Lambda_3(t_{k+1})\cdot (W_{k+1}^1,W_{k+1}^2)
  \right].
\end{equation}

In order to write (\ref{4.24})-(\ref{4.26}) shorter, we introduce the vectors:\\
\quad $\widetilde{\Lambda}_1^1 \in \mathbb{R}^m$ with first components 
$(t_{k+1}-t_{k-1})\Lambda_1^1 (t_{k})$, $1\leq k\leq m-1$
and the last component $(t_{m}-t_{m-1})\Lambda_1^1 (t_{m})$,\\
\quad $\widetilde{\Lambda}_1^2 \in \mathbb{R}^m$ with first components 
$(t_{k+1}-t_{k-1})\Lambda_1^2 (t_{k})$, $1\leq k\leq m-1$
and the last component $(t_{m}-t_{m-1})\Lambda_1^2 (t_{m})$,\\
\quad $\widetilde{\Lambda}_3^1 \in \mathbb{R}^m$ with first components 
$\frac{1}{2}(t_{k+1}-t_{k-1})\Lambda_3^1 (t_{k})$, $1\leq k\leq m-1$
and the last component $\frac{1}{2}(t_{m}-t_{m-1})\Lambda_3^1 (t_{m})$,\\
\quad $\widetilde{\Lambda}_3^2 \in \mathbb{R}^m$ with first components 
$\frac{1}{2}(t_{k+1}-t_{k-1})\Lambda_3^2 (t_{k})$, $1\leq k\leq m-1$
and the last component $\frac{1}{2}(t_{m}-t_{m-1})\Lambda_3^2 (t_{m})$.
Also,  we introduce the vectors in $\mathbb{R}^n$:
\begin{eqnarray*}
\widetilde{\Lambda}_2^1 & = & \frac{1}{2} 
\sum_{0\leq k\leq m-1}(t_{k+1}-t_k)\left( \Lambda_2^1 (t_k)\Phi(Z_k)+\Lambda_2^1 (t_{k+1})\Phi(Z_{k+1})\right)
\\
\widetilde{\Lambda}_2^2 & = & \frac{1}{2}
\sum_{0\leq k\leq m-1}(t_{k+1}-t_k) \left( \Lambda_2^2 (t_k)\Phi(Z_k)+\Lambda_2^2 (t_{k+1})\Phi(Z_{k+1})\right)
\end{eqnarray*}
where $\Phi(Z_k)=(\phi_i(Z_k))_{i\in I}^T \in \mathbb{R}^n$.

\begin{proposition}\label{prop:3.23}
The discrete version of the (\ref{3.23}) is 
\begin{eqnarray}\label{4.26a}
\qquad dJ_{(G,U)}(R,V) & = &  
\left( \left(L(Y)\right)^T M_{ED} +\frac{2}{\epsilon}Y^T N(Z)\right)
  K^{-1}B^1(G,\epsilon)V
\\
&+& \left( \left(L(Y)\right)^T M_{ED} +\frac{2}{\epsilon}Y^T N(Z)\right)
  K^{-1}C^1(G,\epsilon,U)R
\nonumber\\
&+&\frac{1}{\epsilon}\left( (\widetilde{\Lambda}_1^1)^T W^1 + (\widetilde{\Lambda}_1^2)^T W^2 \right)
\nonumber\\
&+&\frac{1}{\epsilon}\left( -(\widetilde{\Lambda}_2^1)^T (\Pi_h^2 R)
+ (\widetilde{\Lambda}_2^2)^T (\Pi_h^1 R)
\right)
\nonumber\\
&+&\frac{1}{\epsilon}\left( (\widetilde{\Lambda}_3^1)^T W^1 + (\widetilde{\Lambda}_3^2)^T W^2 \right).
\nonumber
\end{eqnarray}
\end{proposition}

\noindent
\textbf{Proof.}
We obtain (\ref{4.26a}) by summing (\ref{4.23})-(\ref{4.26}).
\quad$\Box$

\bigskip

Next, we give more details about the relationship between $W$ and $R$.
Let us introduce the $2\times 2$ matrices
$$
M_2(k)=\left(
\begin{array}{cc}
1 -(t_{k+1}-t_k) \partial_1^h\partial_2^h g_h(Z_k) & -(t_{k+1}-t_k) \partial_2^h\partial_2^h g_h(Z_k)\\
(t_{k+1}-t_k) \partial_1^h\partial_1^h g_h(Z_k) & 1 +(t_{k+1}-t_k) \partial_2^h\partial_1^h g_h(Z_k)
\end{array}
\right),
$$
$$
I_2=\left(
\begin{array}{cc}
  1 & 0\\
  0 & 1
\end{array}
\right)  
$$
and the $2\times n$ matrice
$$
N_2(k)=\left(
\begin{array}{r}
-(t_{k+1}-t_k) \Phi^T(Z_k) \Pi_h^2\\
(t_{k+1}-t_k) \Phi^T(Z_k) \Pi_h^1
\end{array}
\right).
$$
We remark that $M_2$ depends on $G$ and $Z$ and $N_2$ on $Z$.
The system (\ref{4.11})-(\ref{4.12}) can be written as 
$$
\left(
\begin{array}{c}
W_{k+1}^1\\
W_{k+1}^2
\end{array}
\right)
=M_2(k)
\left(
\begin{array}{c}
W_{k}^1\\
W_{k}^2
\end{array}
\right)
+N_2(k)R.
$$

\begin{proposition}\label{prop:B23}
We have the following equality
\begin{eqnarray}
\left(
\begin{array}{c}
\begin{array}{c}
W_{1}^1\\
W_{1}^2
\end{array}
\\
\vdots
\\
\begin{array}{c}
W_{m}^1\\
W_{m}^2
\end{array}
\end{array}
\right)
&=& M_{2m}
\times
\left(
\begin{array}{c}
N_2(0)\\
N_2(1)\\
\vdots \\
N_2(m-1)
\end{array}
\right)
R
\label{4.B23}
\end{eqnarray}
where at the right-hand side, $M_{2m}$ is a $2m\times 2m$ matrix defined by
$$
\left(
\begin{array}{ccccc}
I_2 & 0      & \cdots & 0 & 0\\
M_2(1) & I_2 & \cdots & 0& 0\\
\vdots & \vdots& & \vdots& \vdots \\
M_2(m-1)\cdots M_2(1), & M_2(m-1)\cdots M_2(2), & \cdots & M_2(m-1), & I_2
\end{array}
\right)
$$
and the size of the second matrix, which contains $N_2$, is $2m\times n$.
\end{proposition}

\noindent
\textbf{Proof.}
From (\ref{4.13}) and the recurrent relation, we have
\begin{eqnarray*}
\left(
\begin{array}{c}
W_{1}^1\\
W_{1}^2
\end{array}
\right) & = & N_2(0)R\\
\left(
\begin{array}{c}
W_{2}^1\\
W_{2}^2
\end{array}
\right) & = & M_2(1) \left(
\begin{array}{c}
W_{1}^1\\
W_{1}^2
\end{array}
\right)
+N_2(1)R
= M_2(1)N_2(0)R  +  N_2(1)R\\
\vdots &&\\
\left(
\begin{array}{c}
W_{m-1}^1\\
W_{m-1}^2
\end{array}
\right) & = &
M_2(m-2)\cdots M_2(1)N_2(0)R
+ M_2(m-2)\cdots M_2(2)N_2(1)R\\
&&
+\cdots + M_2(m-2)N_2(m-3)R
+ N_2(m-2)R
\\
\left(
\begin{array}{c}
W_{m}^1\\
W_{m}^2
\end{array}
\right) & = &
M_2(m-1)\left(
\begin{array}{c}
W_{m-1}^1\\
W_{m-1}^2
\end{array}
\right)
+N_2(m-1)R\\
&=&
M_2(m-1)M_2(m-2)\cdots M_2(1)N_2(0)R\\
&&
+ M_2(m-1)M_2(m-2)\cdots M_2(2)N_2(1)R \\
&&
+\cdots+ M_2(m-1)M_2(m-2)N_2(m-3)R\\
&&
+ M_2(m-1)N_2(m-2)R + N_2(m-1)R
\end{eqnarray*}
which gives (\ref{4.B23}).
\quad$\Box$

Since $W$ depends on $R$ by (\ref{4.B23}), we can introduce
the linear operator approximation of $B$ in the Corollary \ref{cor:1}
$$
R \in \mathbb{R}^n \rightarrow W=(W^1,W^2)=\left(B^2(G,Z)R,
B^3(G,Z)R \right)\in \mathbb{R}^m \times \mathbb{R}^m.
$$
If we denote by $\ell_i$ the i-th line of the matrix $M_{2m}$ at the right-hand side
of (\ref{4.B23}), for $1\leq i \leq 2m$, then
$$
B^2(G,Z)=
\left(
\begin{array}{c}
\ell_1\\
\ell_3\\
\vdots \\
\ell_{2m-1}
\end{array}
\right)
\left(
\begin{array}{c}
N_2(0)\\
N_2(1)\\
\vdots \\
N_2(m-1)
\end{array}
\right),
$$
$$
B^3(G,Z)=
\left(
\begin{array}{c}
\ell_2\\
\ell_4\\
\vdots \\
\ell_{2m}
\end{array}
\right)
\left(
\begin{array}{c}
N_2(0)\\
N_2(1)\\
\vdots \\
N_2(m-1)
\end{array}
\right)
$$
and $B^2(G,Z)$, $B^3(G,Z)$ are $m\times n$ matrices.
The size of the matrix containing $N_2$ is $2m\times n$.

\subsection{Gradient type algorithm}

We start by presenting the algorithm.

\bigskip

\textbf{Step 1} Start with $k=0$, $\epsilon >0$ some given ``small'' parameter and select
some initial $(G^k,U^k)$.

\textbf{Step 2} Compute $Y^k$ the solution of (\ref{4.2}) and $Z^k$ solution of
(\ref{4.3})-(\ref{4.5}).

\textbf{Step 3} Find $(R^k,V^k)$ such that $dJ_{(G^k,U^k)}(R^k,V^k) <0$.
We say that $(R^k,V^k)$ is a descent direction.

\textbf{Step 4} Define $(G^{k+1},U^{k+1})=(G^k,U^k)+\lambda_k (R^k,V^k)$,
where $\lambda_k >0$ is obtained via some line search
$$
\lambda_k \in \arg\min_{\lambda >0} J\left((G^k,U^k)+\lambda (R^k,V^k)\right).
$$

\textbf{Step 5} If $| J(G^{k+1},U^{k+1}) - J(G^k,U^k)|$ is  below some prescribed 
tolerance parameter, then \textbf{Stop}.
If not, update $k:=k+1$ and go to \textbf{Step 3}.

\bigskip

In the \textbf{Step 3}, we have to provide a descent direction.

We present in the following a partial result.

Let us introduce a simplified adjoint system: find $p_h\in \mathbb{V}_h$ such that
\begin{equation}\label{4.27}
  \int_D \nabla \varphi_h \cdot \nabla p_h d\mathbf{x}
  = \int_E \partial_2 j(\mathbf{x}, y_h(\mathbf{x}))\varphi_h d\mathbf{x}
  +\frac{2}{\epsilon}\int_0^{T_g} y_h(Z(t))\varphi_h(Z(t))|Z^\prime(t)| dt
\end{equation}
for all $\varphi_h \in \mathbb{V}_h$ and with $Z(t)$ given by (\ref{4.3})-(\ref{4.5}).
We have $p_h(\mathbf{x})=\sum_{i\in I_0} P_i \phi_i(\mathbf{x})$ and
$P=(P_i)_{i\in I_0}^T\in \mathbb{R}^{n_0}$.
The linear system associated to (\ref{4.27}) is
$$
KP=M_{ED}^T L(Y)+\frac{2}{\epsilon} N(Z)Y.
$$
We recall that $K$ and $N(Z)$ are symmetric matrices.

\begin{proposition}\label{prop:4.1}
  Given $g_h,u_h\in\mathbb{W}_h$, let $y_h\in\mathbb{V}_h$ the solution of
  (\ref{4.1}). If $r_h=-p_hu_h$ and $v_h=-p_h$, where $p_h\in\mathbb{V}_h$
  is the solution of
  (\ref{4.27}), then
\begin{equation}\label{4.28}
\int_E \partial_2 j(\mathbf{x}, y_h(\mathbf{x})) q_h d\mathbf{x}
+\frac{2}{\epsilon}\int_0^{T_g} y_h(Z(t))q_h(Z(t))|Z^\prime(t)| dt
\leq 0,
\end{equation}
where $q_h\in\mathbb{V}_h$ is the solution of (\ref{4.7}) depending on $r_h$
and $v_h$.
\end{proposition}  

\noindent
\textbf{Proof.}
Putting $\varphi_h=p_h$ in (\ref{4.7}) and $\varphi_h=q_h$ in (\ref{4.27}), we get
\begin{eqnarray}\label{4.29}
  && \int_D \left( (g_h+\epsilon)_+^2 v_h + 2(g_h+\epsilon)_+ u_h r_h\right)
  p_h d\mathbf{x}
  =\int_D \nabla q_h \cdot \nabla p_h d\mathbf{x}
  \\
  &=&\int_E \partial_2 j(\mathbf{x}, y_h(\mathbf{x}))q_h d\mathbf{x}
  +\frac{2}{\epsilon}\int_0^{T_g} y_h(Z(t))q_h(Z(t))|Z^\prime(t)| dt .
  \nonumber
\end{eqnarray}

For $v_h=-p_h$, we have
$$
\int_D (g_h+\epsilon)_+^2 v_hp_hd\mathbf{x}
=-\int_D (g_h+\epsilon)_+^2 p_h^2 d\mathbf{x}
\leq 0
$$
and for $r_h=-p_hu_h$, we have
$$
\int_D  2(g_h+\epsilon)_+ u_h r_h p_h d\mathbf{x}
=-\int_D  2(g_h+\epsilon)_+ (u_h p_h)^2 d\mathbf{x}
\leq 0
$$
since $(g_h+\epsilon)_+ \geq 0$ in $D$.
This ends the proof.
\quad$\Box$

\begin{remark}\label{rem:4.2}
The left-hand side of (\ref{4.28}) represents the first two terms of
  (\ref{4.22}). 
We can obtain a similar result as in Proposition \ref{prop:4.1}, without using the adjoint system,
by taking
\begin{eqnarray*}
(V^*)^T &=& - \left( \left(L(Y)\right)^T M_{ED} +\frac{2}{\epsilon}Y^T N(Z)\right)
  K^{-1}B^1(G,\epsilon) \\
(R^*)^T &=& - \left( \left(L(Y)\right)^T M_{ED} +\frac{2}{\epsilon}Y^T N(Z)\right)
  K^{-1}C^1(G,\epsilon,U),
\end{eqnarray*}
in place of $V$ and $R$ in (\ref{4.23}). In this case, (\ref{4.23}) becames
$-\| V^*\|_{\mathbb{R}^n}^2 -\| R^*\|_{\mathbb{R}^n}^2 \leq 0$.
We point out that $(V^*)^T=-P^TB^1(G,\epsilon)$ and
$(R^*)^T=-P^TC^1(G,\epsilon,U)$, so $(V^*,R^*)$ is different from the direction given by
Proposition \ref{prop:4.1}.
\end{remark}

Now, we present a descent direction, obtained from the complete gradient of the discrete cost (\ref{4.22}).

\begin{proposition}\label{prop:4.2}
For $(R^{**},V^*)\in \mathbb{R}^n \times \mathbb{R}^n$ given by
\begin{eqnarray*}
(V^*)^T&=& - \left( \left(L(Y)\right)^T M_{ED} +\frac{2}{\epsilon}Y^T N(Z)\right)
  K^{-1}B^1(G,\epsilon) \\
(R^{**})^T&=& - \left( \left(L(Y)\right)^T M_{ED} +\frac{2}{\epsilon}Y^T N(Z)\right)
  K^{-1}C^1(G,\epsilon,U) \\
&& - \frac{1}{\epsilon}\left( (\widetilde{\Lambda}_1^1)^T B^2(G,Z) + (\widetilde{\Lambda}_1^2)^T B^3(G,Z) \right)
\nonumber\\
&& -\frac{1}{\epsilon}\left( -(\widetilde{\Lambda}_2^1)^T \Pi_h^2
+ (\widetilde{\Lambda}_2^2)^T \Pi_h^1 
\right)
\nonumber\\
&& -\frac{1}{\epsilon}\left( (\widetilde{\Lambda}_3^1)^T B^2(G,Z) + (\widetilde{\Lambda}_3^2)^T B^3(G,Z) \right).
\nonumber
\end{eqnarray*}
we obtain a descent direction for $J$ at $(G,U)$.
\end{proposition}  

\noindent
\textbf{Proof.}
In (\ref{4.26a}), we replace $W^1$ by $B^2(G,Z)R$ and $W^2$ by $B^3(G,Z)R$, we obtain that
$dJ_{(G,U)}(R,V)=-(V^*)^T V - (R^{**})^T R$, then
$dJ_{(G,U)}(R^{**},V^*)=-\| V^*\|_{\mathbb{R}^n}^2 -\| R^{**}\|_{\mathbb{R}^n}^2 \leq 0$.
\quad$\Box$

\section{Numerical tests}
\setcounter{equation}{0}

Shape optimization problems and their penalization are strongly nonconvex. The computed optimal domain depends on the starting domain, but also on the penalization
$\epsilon$ or other numerical parameters. It may be just a local optimal solution. 

Moreover, the final computed value of the penalization integral is small, but not null. This allows differences between the  optimal  computed domain $\Omega_g$ and the zero level curves of the computed optimal state $y_{\epsilon}$. Consequently, we compare the obtained optimal cost in the penalized problem with the costs  in the original problem (\ref{1.1}) - (\ref{1.3}) corresponding to the optimal  computed domain $\Omega_g$ and the zero level curves of $y_{\epsilon}$. This is a standard procedure, to inject the approximating optimal solution in the original problem. Notice that in all the experiments, the cost corresponding to $\Omega_g$ is the best one, but the differences with respect to the other computed cost values are small. This shows that the rather complex approximation/penalization that we use is reasonable. Its advantage is that it may be used as well in the case of boundary observation or for Neumann boundary conditions and this will be performed in a subsequent paper.

In the examples, we have employed the software FreeFem++, \cite{freefem++}.
 
\medskip

\textbf{Example 1.}

The computational domain is $D=]-3,3[\times ]-3,3[$ and
the observation zone $E$ is the disk of center $(0,0)$ and radius $0.5$.
The load is $f=1$, $j(g)=(y_\epsilon-y_d)^2$, where
$y_d(x_1,x_2)=-(x_1-0.5)^2 -(x_2-0.5)^2+\frac{1}{16}$, then the cost function (\ref{3.2}) becomes
\begin{equation}\label{5.1}
\min_{g\in \mathcal{F},\, u\in L^2(D)} J(g,u)=
\left\{ 
\int_E (y_\epsilon-y_d)^2d\mathbf{x}
+\frac{1}{\epsilon}\int_{I_g} \left(y_\epsilon(\mathbf{z}_g(t))\right)^2|\mathbf{z}_g^\prime(t)|dt
\right\} .
\end{equation}

The mesh of $D$ has 73786 triangles and 37254 vertices.
The penalization parameter is $\epsilon=10^{-3}$ and the tolerance
parameter for the stopping test
at the \textbf{Step 5} of the algorithm is $tol=10^{-6}$.
The initial domain is the disk of center $(0,0)$ and radius $2.5$
with a circular hole of center $(-1,-1)$ and radius $0.5$.

At the \textbf{Step 3} of the Algorithm,
we use $(R^k,V^k)$ given by Proposition \ref{prop:4.1}.
At the \textbf{Step 4}, in order to have $E\subset \Omega_k$,
we use a projection $\mathcal{P}$ at the line search
$$
\lambda_k \in \arg\min_{\lambda >0}
J\left( \mathcal{P}(G^k+\lambda R^k),
U^k+\lambda V^k\right)
$$
and $G^{k+1}=\mathcal{P}(G^k+\lambda_k R^k)$.
If the value of $g_h^{k}+\lambda r_h^k$ at a vertex from $\overline{E}$
is positive, then we set this value to $-0.1$.
We recall that the left-hand side of (\ref{4.28}) represents only the first two
terms of (\ref{4.22}), not the whole gradient.
\begin{figure}[ht]
\begin{center}
\includegraphics[width=5.5cm]{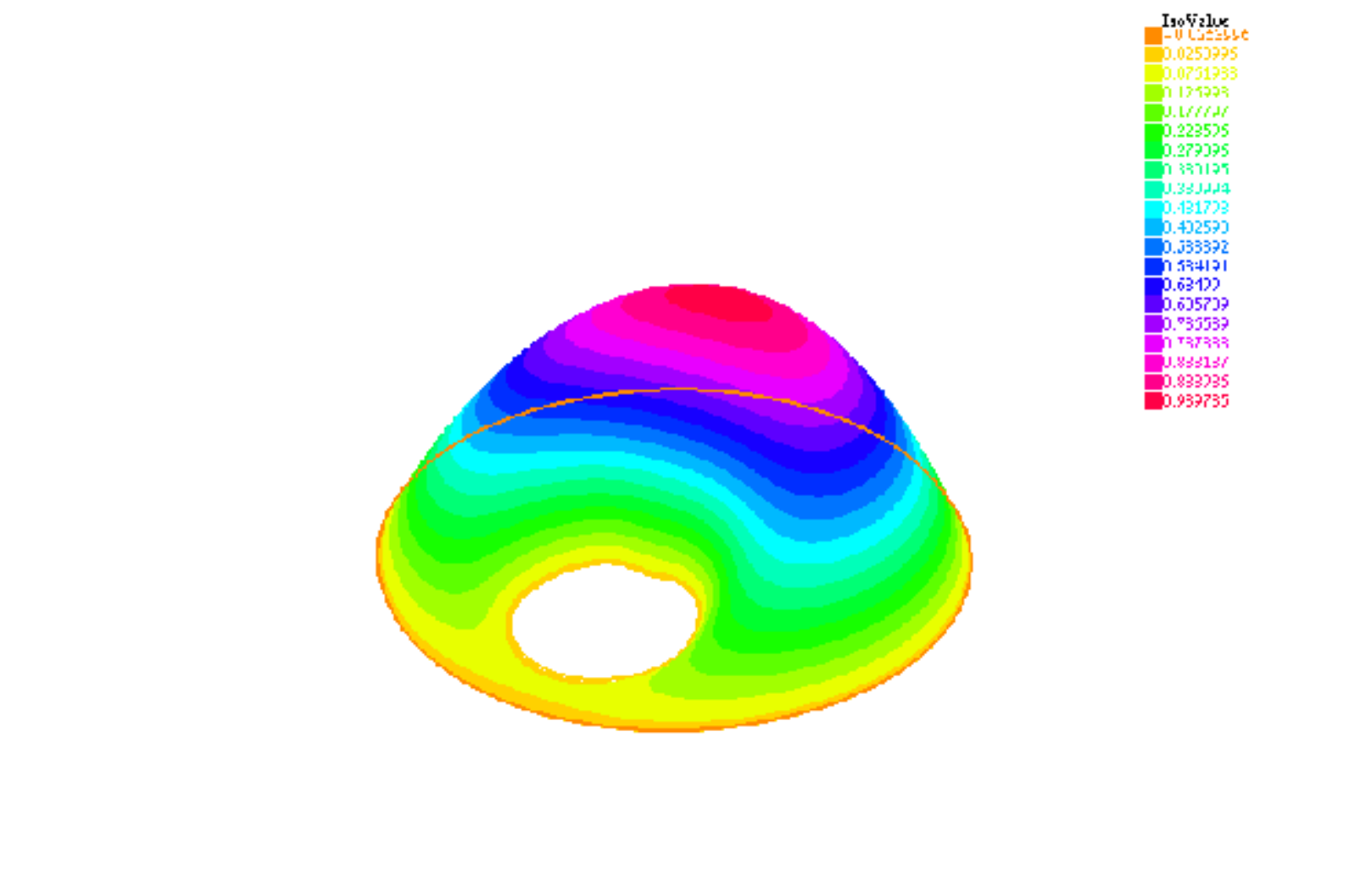}
\
\includegraphics[width=5.5cm]{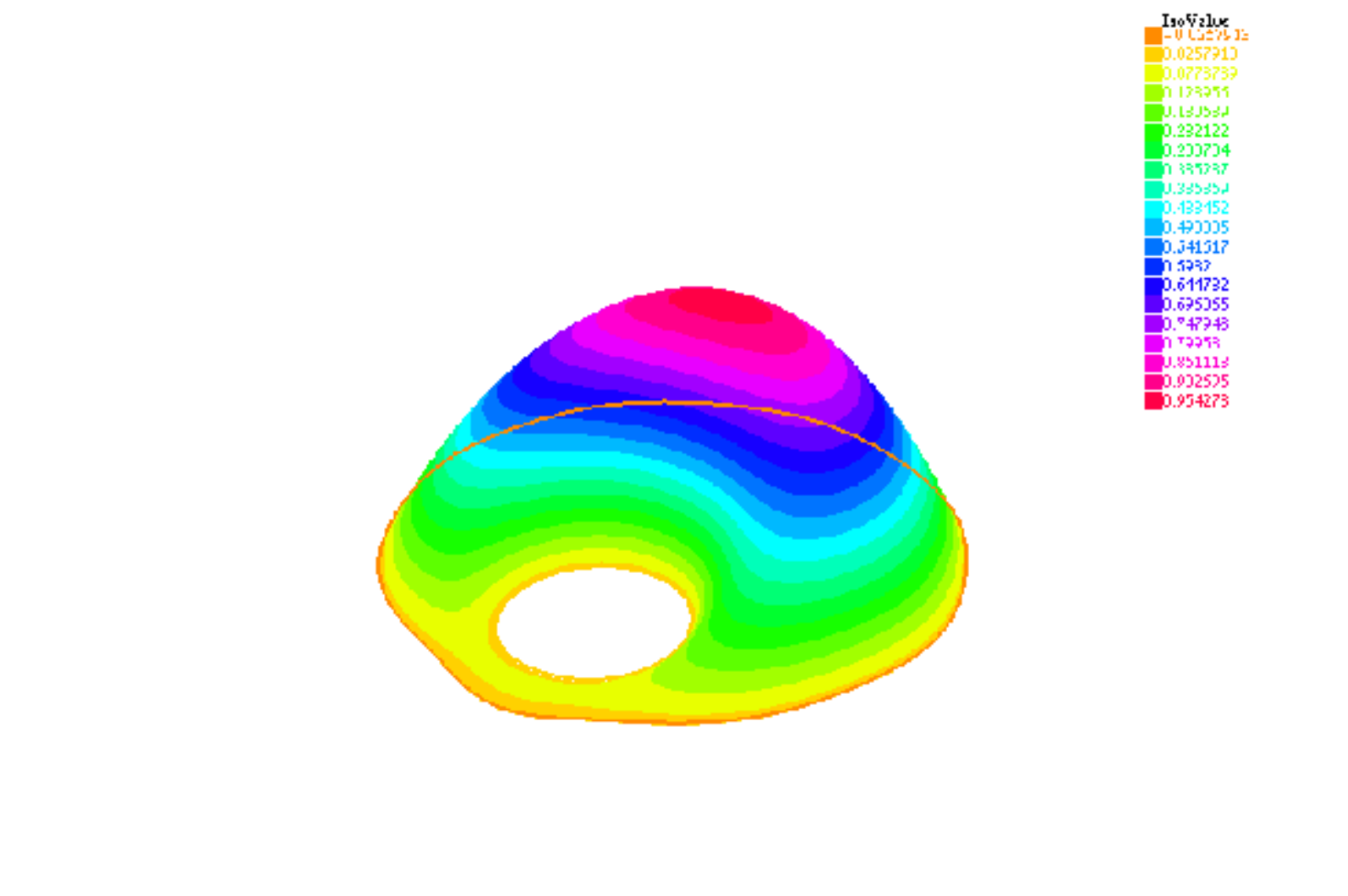}
\
\includegraphics[width=5.5cm]{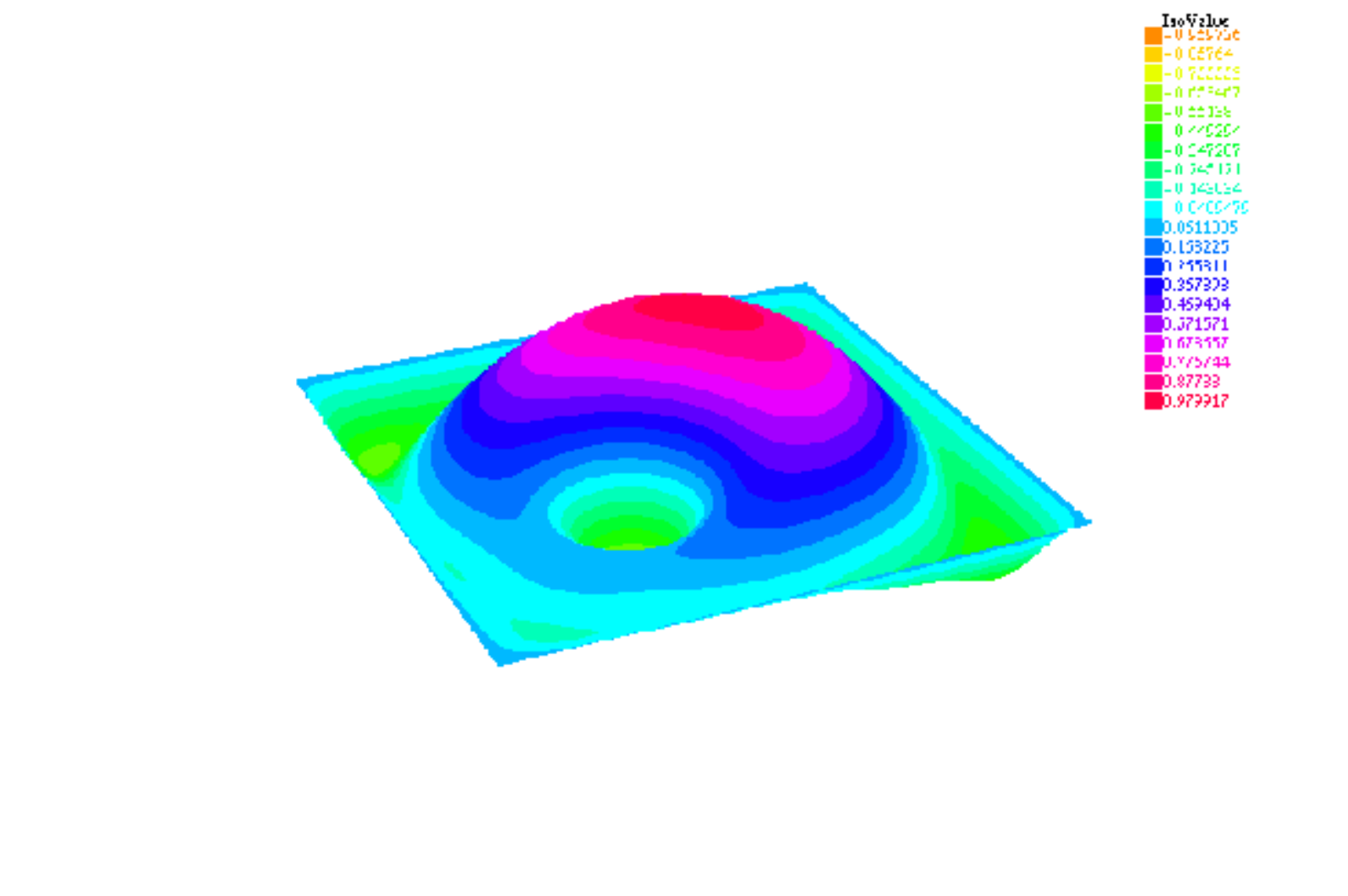}
\end{center}
\caption{Example 1. 
The solution of the elliptic problem (\ref{1.2})-(\ref{1.3})
in the domain $\Omega_g$ (left),
in the domain bounded by the zero level sets of $y_\epsilon$ (right) and
the final computed state $y_\epsilon$ in $D$ (bottom).
\label{fig:ex2_y}}
\end{figure}
If $r_h,\ v_h$ are given by Proposition \ref{prop:4.1} and $\gamma >0$ is a scaling
parameter, then $\gamma r_h$ and $v_h$ verify (\ref{4.28}), that is they also give a descent direction.
We take the scaling parameter for $r_h$ given by $\gamma=\frac{1}{\max(r_h)}$, that is a normalization of $r_h$. 
In this way we avoid the appearance of very high values of the objective function, that may stop the 
algorithm even in the first iteration.  For the line search at the \textbf{Step 4}, we use 
$\lambda=\rho^i\lambda_0$, with $\lambda_0=1$, $\rho=0.5$ for $i=0,1,\dots,30$.

The stopping test is obtained for $k=94$ and some values of the objective
function are:
$J(G^0,U^0)=33110.5$, $J(G^{30},U^{30})=54.725$, $J(G^{94},U^{94})=14.9851$.

At the final iteration, the first term of the
optimal objective function is $1.03796$ and
$\int_{\partial\Omega_g} y_\epsilon^2(s)ds=1.39471 \times 10^{-2}$. 
We point out that the optimal $\Omega_g$ has a hole and the penalization term is a sum of two integrals
$$
\int_{\partial\Omega_g} y_\epsilon^2(s)ds=
\sum_{j=1}^2\int_{I_j} \left(y_\epsilon(\mathbf{z}_g(t))\right)^2|\mathbf{z}_g^\prime(t)|dt
$$
where the integral over $I_1$ corresponding to the exterior boundary of $\Omega_g$
and $I_2$ to the boundary of the hole.
In Figure \ref{fig:ex2_y} in the bottom, we can see the computed optimal state $y_\epsilon$ in iteration 94.
We also compute the costs $\int_E (y_1-y_d)^2d\mathbf{x}=0.998189$ where $y_1$ is the solution
of the initial elliptic problem (\ref{1.2})-(\ref{1.3})
in the domain $\Omega_g$ with g obtained in iteration 94
and $\int_E (y_2-y_d)^2d\mathbf{x}=1.04032$ where $y_2$ is the solution
of the elliptic problem (\ref{1.2})-(\ref{1.3})
in the domain bounded by the zero level sets of $y_\epsilon$ in iteration 94,
see Figure~\ref{fig:ex2_y}.

\medskip
\textbf{Example 2.}

The domains $D$, $E$ and the mesh of $D$ are the same as in Example 1.
For $f=4$ and $y_d(x_1,x_2)=-x_1^2 -x_2^2+1$, 
we have the exact optimal state $y=y_d$ defined in the disk of
center $(0,0)$ and radius $1$, that gives an optimal domain of the problem (\ref{1.1})-(\ref{1.3}).

We have used $\epsilon=10^{-1}$ and the starting configuration:
the disk of center $(0,0)$ and radius $2.5$
with the circular hole of center $(-1,-1)$ and radius $0.5$.
We use $(R^k,V^k)$ given by Proposition \ref{prop:4.1}.
The parameters for the line search and $\gamma$ are the same as in the precedent example.

The stopping test is obtained for $k=64$.
The initial and the final computed values of the objective function
are $5368.84$ and $11.2311$.
We obtain a local minimum that is different from the above global solution. The first term of the
final computed objective function is $0.472856$.
The term $\int_{\partial\Omega_g} y_\epsilon^2(s)ds$ is $1.07583$ and it was computed over the exterior
boundary as well as over the boundaries of two holes. The length of the total boundary
of the optimal domain is
$23.9714$ and of the initial domain is $2\pi(2.5+0.5)=18.8495$.

\begin{figure}[ht]
\begin{center}  
\includegraphics[width=5.5cm]{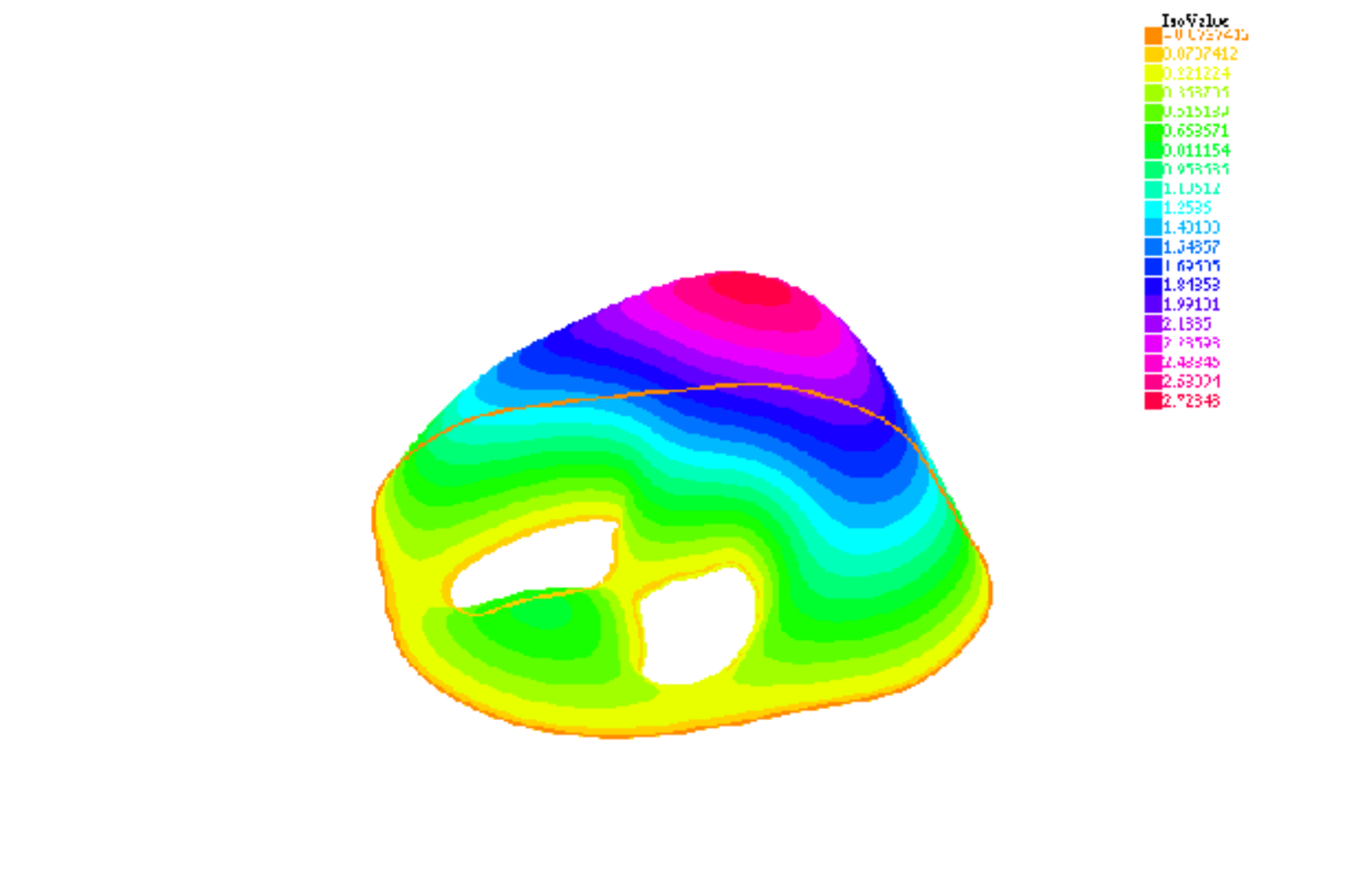}
\
\includegraphics[width=5.5cm]{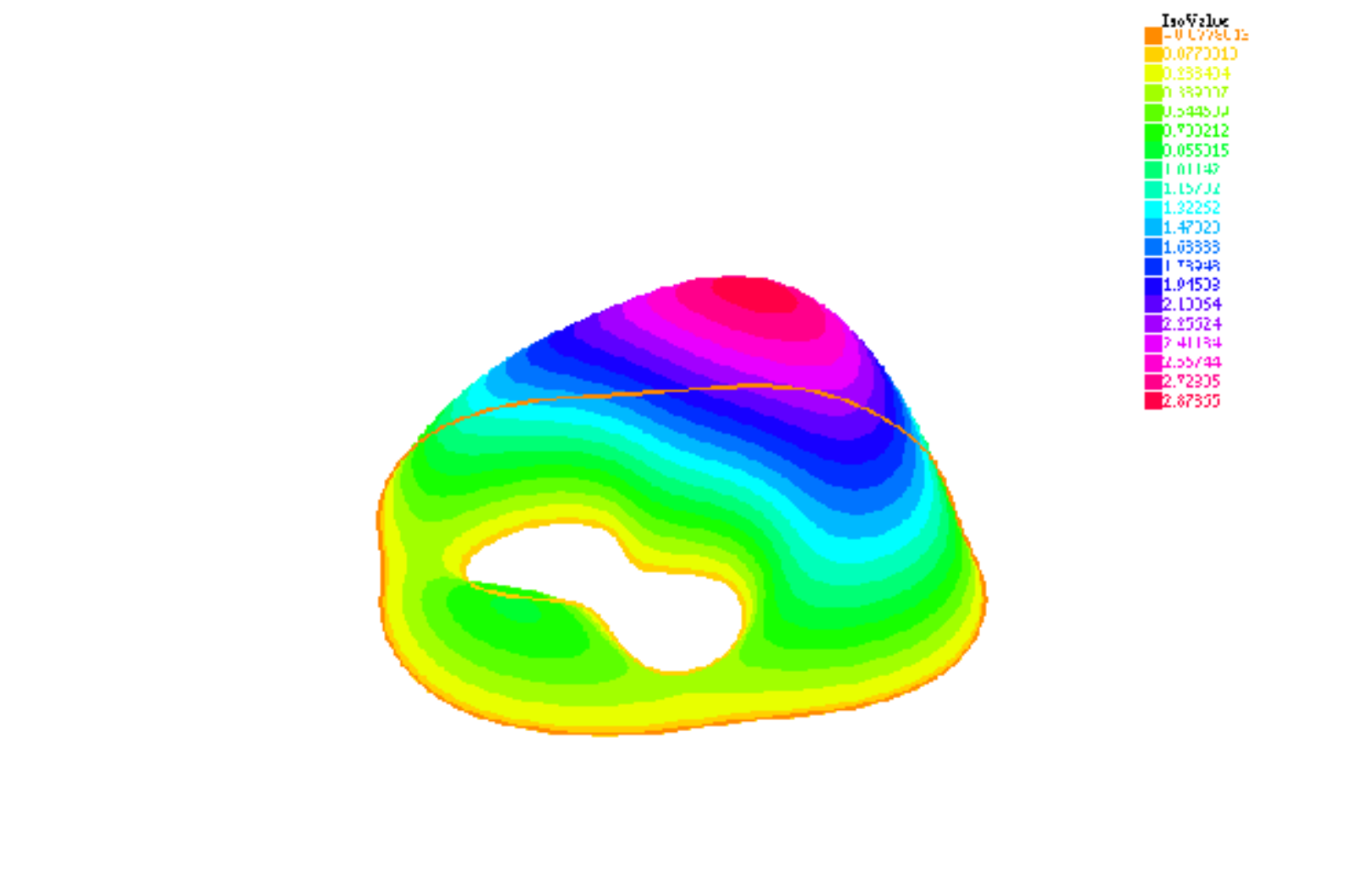}
\
\includegraphics[width=5.5cm]{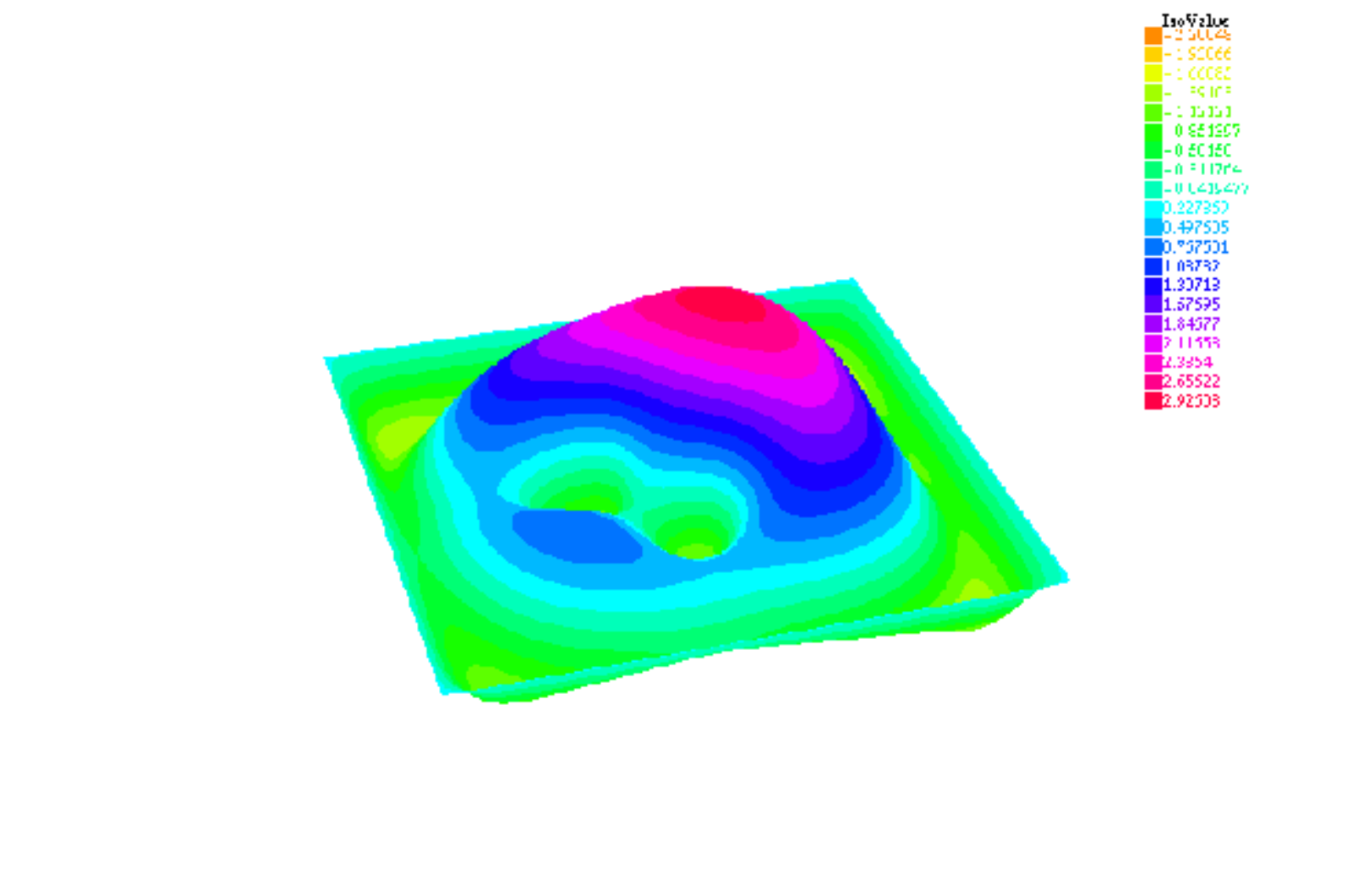}
\end{center}
\caption{Example 2. 
The numerical solution of the elliptic problem (\ref{1.2})-(\ref{1.3})
in the optimal domain $\Omega_g$ (left),
in the domain bounded by the zero level sets of $y_\epsilon$ (right) and
the computed optimal state $y_\epsilon$ (bottom).
\label{fig:ex3_1hole_y_exact}}
\end{figure}

The domain changes its topology.
The computed optimal state $y_\epsilon$ is presented in Figure \ref{fig:ex3_1hole_y_exact} in the bottom.
At the left, we show $y_1$ the solution of the elliptic problem (\ref{1.2})-(\ref{1.3})
in the domain $\Omega_g$ which gives $\int_E (y_1-y_d)^2d\mathbf{x}=0.295178$,
at the right we show $y_2$ the solution of the elliptic problem (\ref{1.2})-(\ref{1.3})
in the domain bounded by the zero level sets of $y_\epsilon$,
which gives $\int_E (y_2-y_d)^2d\mathbf{x}=0.471788$.

\medskip
\textbf{Example 3.}

We have also used the descent direction given by Proposition \ref{prop:4.2},
for the starting configuration the disk of center $(0,0)$ and radius $1.5$,
$\epsilon=10^{-1}$, $\gamma=\frac{1}{\|r_h\|_\infty}$ and a mesh of $D$ of 32446 triangles
and 16464 vertices.
For solving the ODE systems
(\ref{4.3})-(\ref{4.5}) and (\ref{4.11})-(\ref{4.13}) we use $m=30$.

\begin{figure}[ht]
\begin{center}
\includegraphics[width=4.5cm]{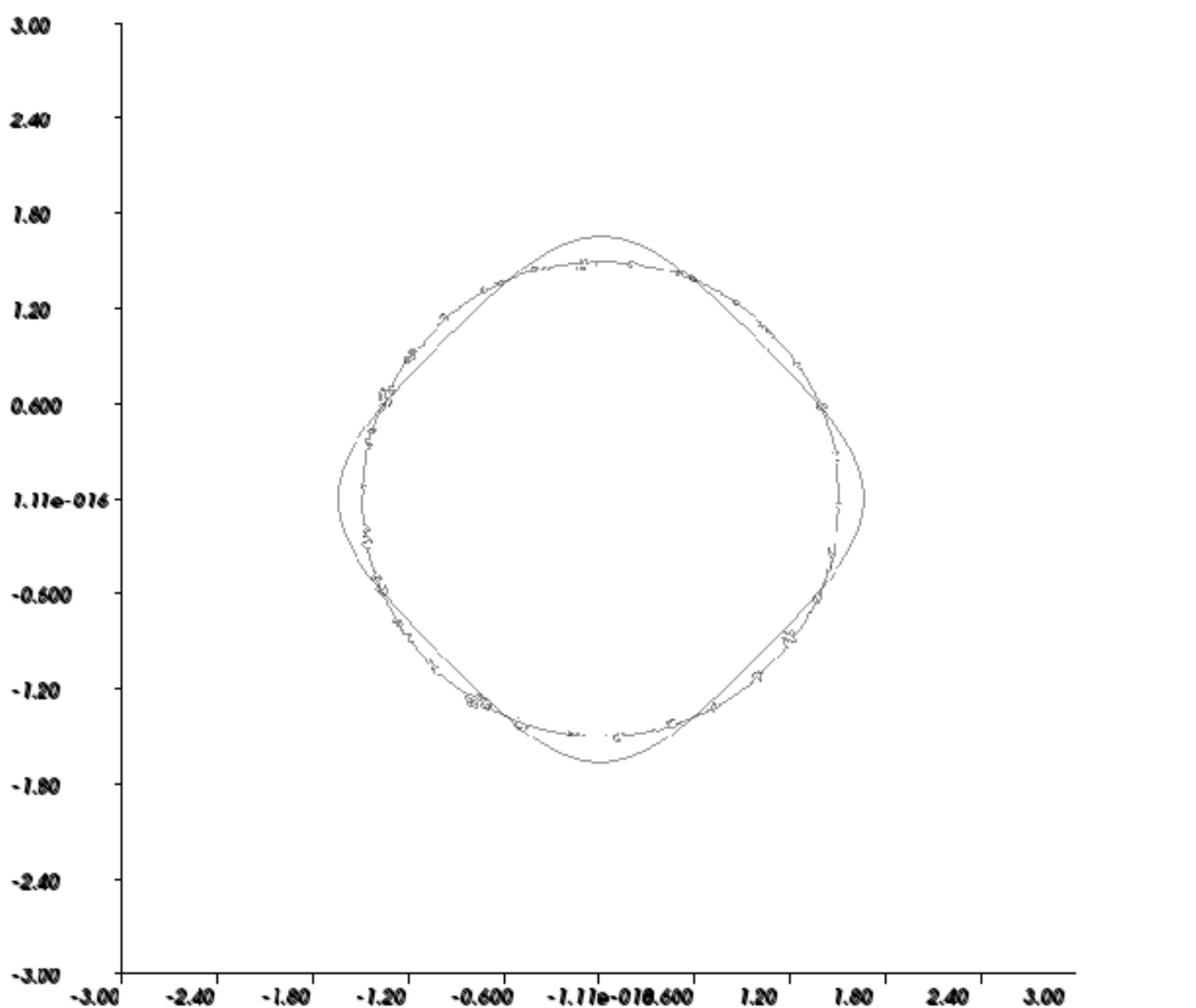}
\   
\includegraphics[width=6.5cm]{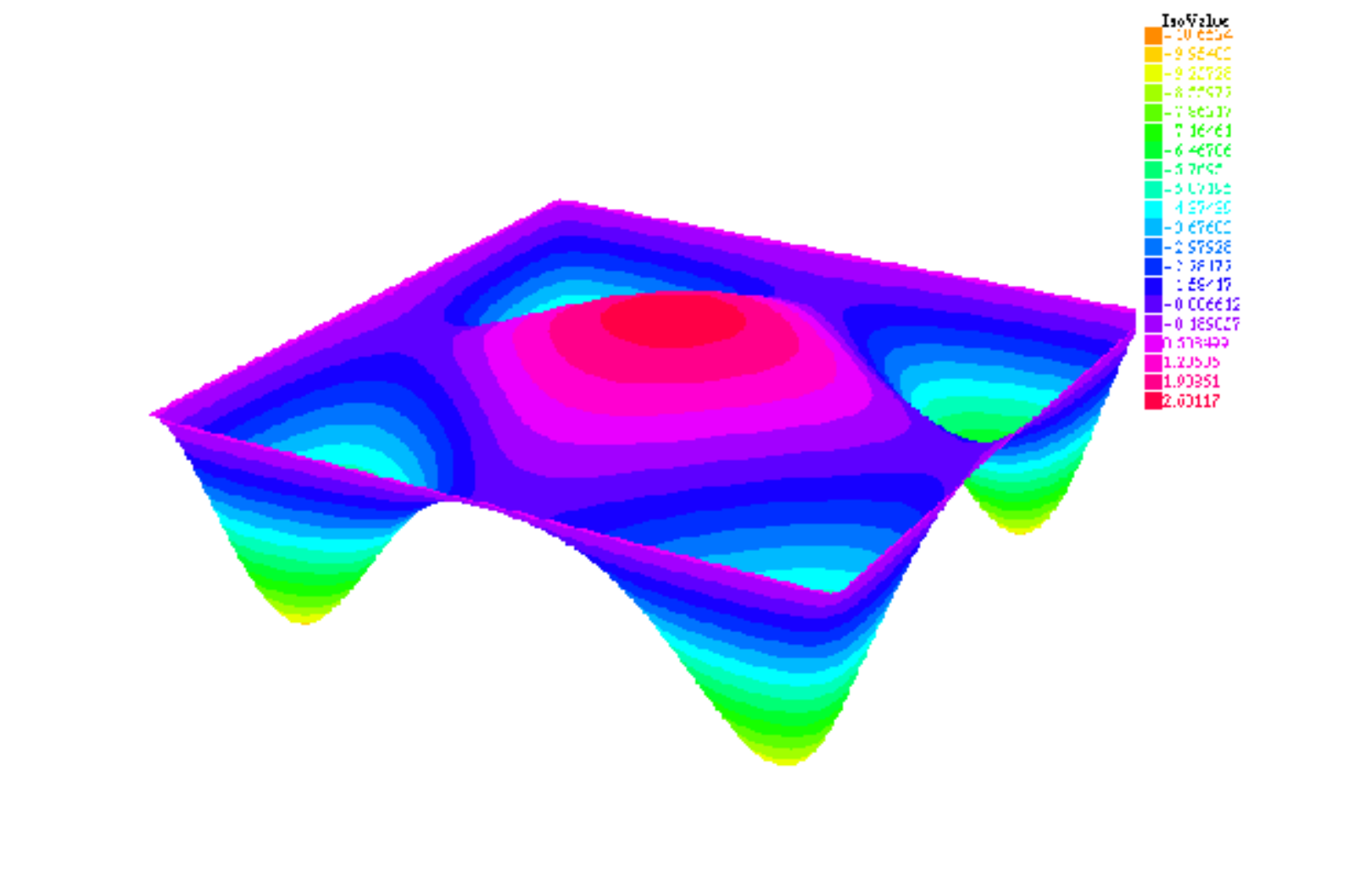}
\
\includegraphics[width=5.5cm]{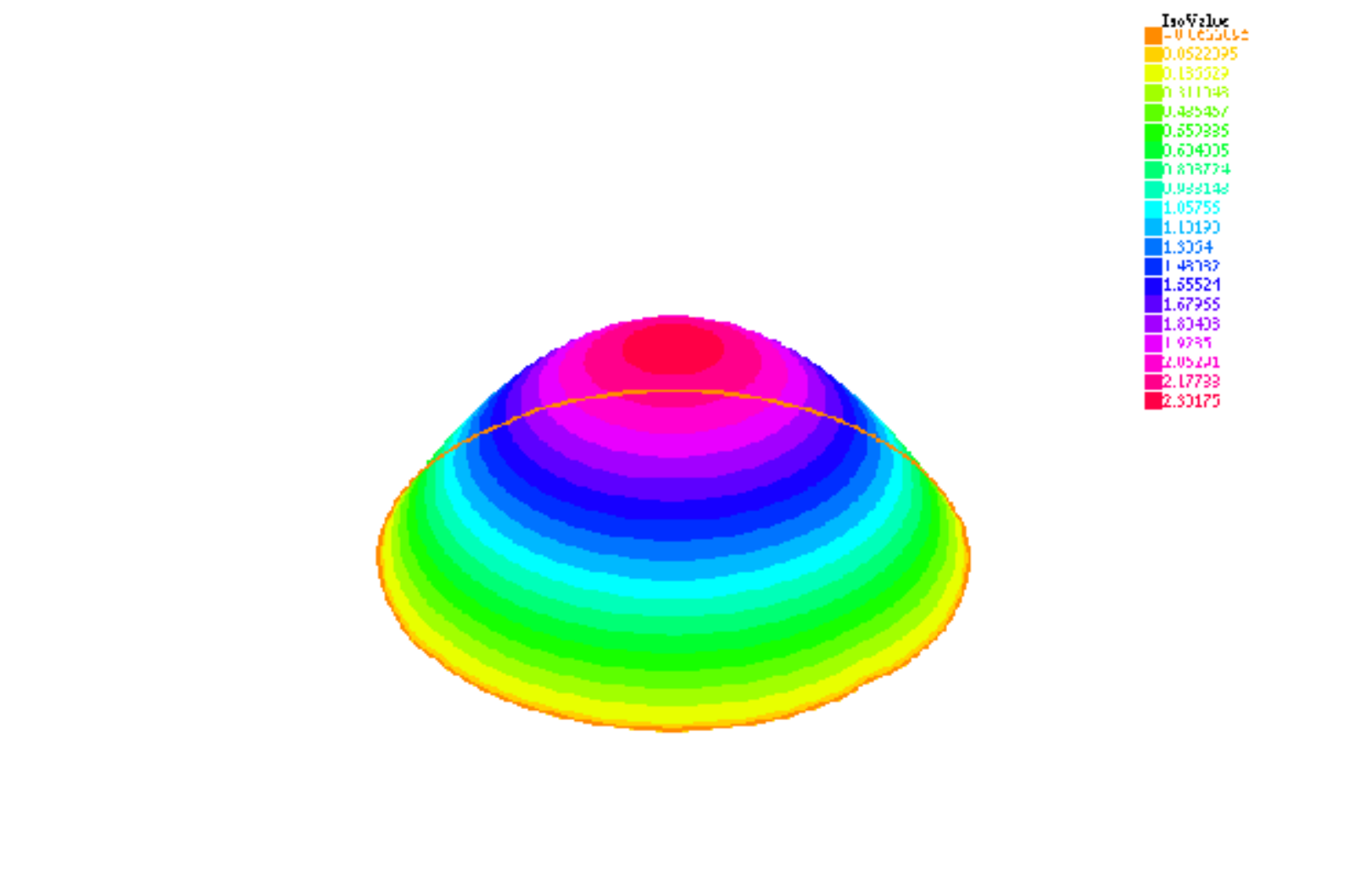}
\
\includegraphics[width=5.5cm]{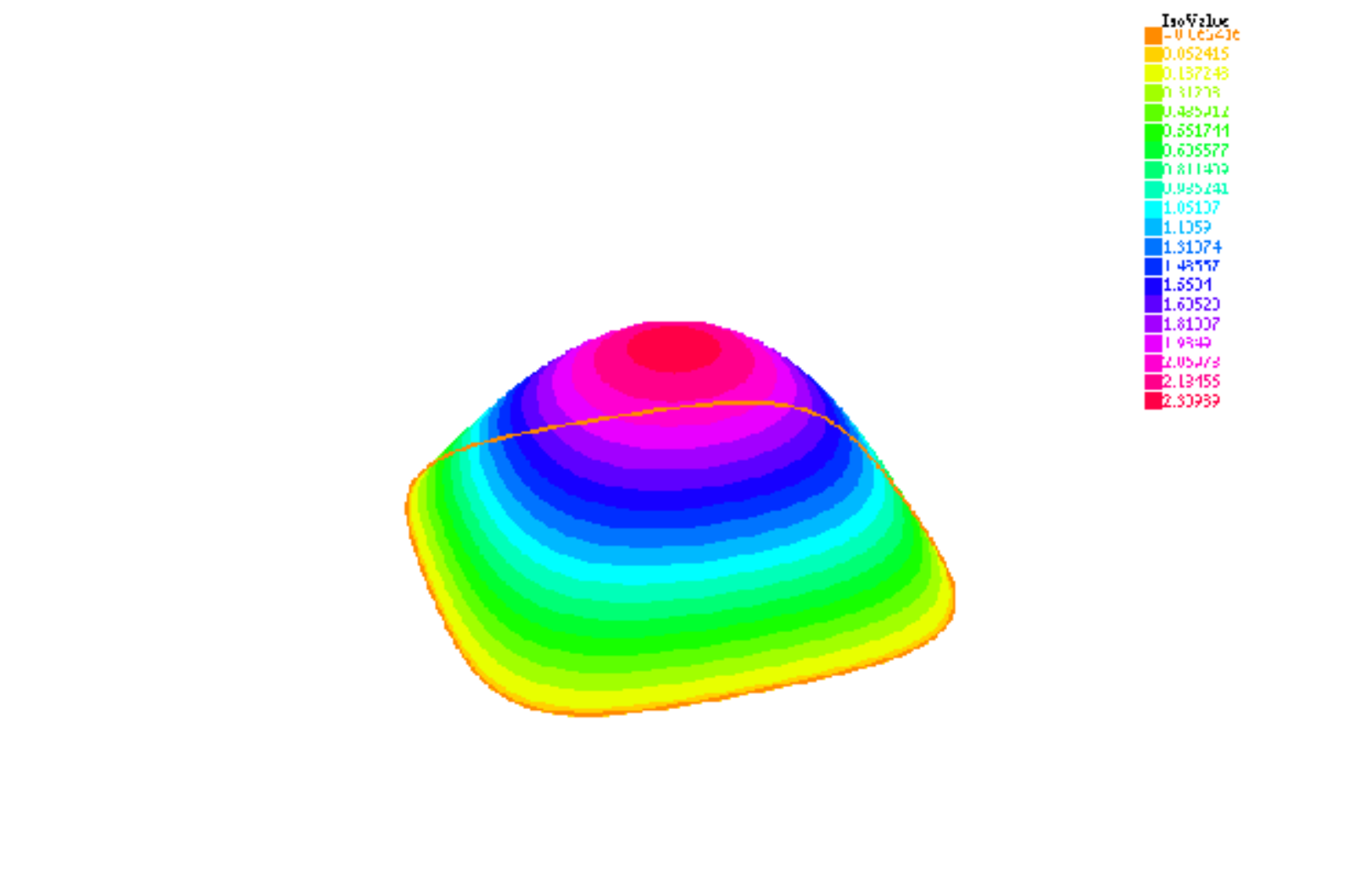}
\end{center}
\caption{Example 3. The zero level sets of the computed optimal $g$, $y_\epsilon$ (top, left),
the final state $y_\epsilon$ (top, right),
the solution of the elliptic problem (\ref{1.2})-(\ref{1.3})
in the domain $\Omega_g$ (bottom, left) and
in the domain bounded by the zero level sets of $y_\epsilon$ (bottom, right).
\label{fig:total_k_20_bnd_g_y_bis}}
\end{figure}

At the initial iteration, we have 
$\int_E (y_\epsilon-y_d)^2d\mathbf{x}=72.3767$, 
$\int_{\partial\Omega_g} y_\epsilon^2(s)ds=658.459$ and the value of the objective function is
$J_0=6656.98$. The algorithm stops after 12 iterations and we have at the final iteration
$\int_E (y_\epsilon-y_d)^2d\mathbf{x}=1.22861$,
$\int_{\partial\Omega_g} y_\epsilon^2(s)ds=0.557556$ and the value of the penalized objective function is
$J_{12}=6.80521$. The final domain is a perturbation of the initial one,
the circular non-smooth curve in the
top, left image of Figure \ref{fig:total_k_20_bnd_g_y_bis}.
We have $\int_E (y_1-y_d)^2d\mathbf{x}=1.20398$ for
$y_1$ the solution of the elliptic problem (\ref{1.2})-(\ref{1.3})
in the final domain $\Omega_g$ and 
$\int_E (y_2-y_d)^2d\mathbf{x}=1.21767$ for 
$y_2$ the solution of the elliptic problem (\ref{1.2})-(\ref{1.3})
in the domain bounded by the zero level sets of $y_\epsilon$,
Figure \ref{fig:total_k_20_bnd_g_y_bis}
at the bottom, right.

Finally, we notice that the hypothesis of Proposition \ref{prop:3.3} is obviously
fulfilled  by the null level sets of $y_\epsilon$ with a corresponding parametrization. In an approximate sense, it is also
fulfilled by the computed optimal domain $\Omega_g$ since the penalization
integral is "small" in all the examples. 

\newpage

\end{document}